\documentclass[twoside,11pt]{article}
\usepackage{graphicx}
\usepackage{amssymb}
\usepackage{amsthm}
\newfont{\bbb}{msbm10 scaled\magstep 1}
\newcommand{\bd}{{\rm bd}}
\newcommand{\inter}{{\rm int}}
\newcommand{\diam}{{\rm diam}}
\newcommand{\conv}{{\rm conv}}
\newcommand{\width}{{\rm width}}

\newtheorem{thm}{Theorem}
\newtheorem{cor}{Corollary}
\newtheorem{pro}{Proposition}
\newtheorem{cla}{Claim}
\newtheorem{lem}{Lemma}

\usepackage{graphicx}
\usepackage{amssymb}
\usepackage{amsthm}
\usepackage{graphicx}
\usepackage{amssymb}
\usepackage{amsthm}

\theoremstyle{remark}
\newtheorem{exa}{\bf Example} 

\newtheorem{prob}{\bf Problem}

\date{}

\title
{\bf Spherical geometry - a survey on width and thickness of convex bodies}

\date{}

\begin{document}

\maketitle
\thispagestyle{empty}

\vskip-1cm
\centerline
{MAREK LASSAK}
\pagestyle{myheadings} \markboth {\centerline {Marek Lassak}}{\centerline {Width and thickness of spherical convex bodies - a survey}}

\baselineskip17.65pt        

\maketitle

\vskip0.9cm
\begin{abstract}
\noindent
This chapter concerns the geometry of convex bodies on the $d$-dimensional sphere $S^d$.  
We concentrate on the results based on the notion of width of a convex body $C \subset S^d$ determined by a supporting hemisphere of $C$.
Important tools are the lunes containing $C$.
The supporting hemispheres take over the role of the supporting half-spaces of a convex body in Euclidean space, and lunes  the role of strips.
Also essential is the notion of thickness of $C$, i.e., its minimum width. 
In particular, we describe properties of reduced spherical convex bodies and spherical bodies of constant width. 
The last notion coincides with the notions of complete bodies and bodies of constant diameter on $S^d$.
The results reminded and commented on here concern mostly the width, thickness, diameter, perimeter, area and extreme points of spherical convex bodies, reduced bodies and bodies of constant width. 
\end{abstract}



\vskip0.2cm
\noindent
The final version of this preprint will appear in the book {\it Surveys in Geometry} (ed. A. Papadopoulos), Springer, 2021.

{\ }

\vskip-0.5cm

\section{Introduction}

\vskip0.16cm
\noindent
Spherical geometry has its sources in the research on spherical trigonometry which in turn was motivated mostly by navigation needs and astronomic observations.
This research started in ancient times mostly by Greeks and continued in the Islamic World,
and also by ancient Chinese, Egyptians and Indians.
The next important stage began in the XVIII-th century by mostly Euler, and by his collaborators and followers. 
Descriptions of these achievements is presented by Papadopoulos \cite{Pap2}. 
A historical look to spherical geometry is given by Van Brummelen in his half-and-half historical and popular science book \cite{VB}, by Rosenfeld \cite{Ro} and Whittlesey \cite{Wh}.

There are many textbooks on spherical trigonometry, as for instance those from the XIX-th century by Todhunter \cite{Tod} and Murray \cite{Mu}.
The oldest important research works on spherical convex bodies that the author of this chapter found are those of Kubota \cite{Kub} and Blaschke~\cite{Bla1}.

The aim of this chapter is to summarize and comment on properties of spherical convex bodies, mostly these properties which are consequences of the notion of the width of a spherical convex body $C$ determined by a supporting hemisphere of $C$.
This notion of width is introduced in \cite{L3}.
We omit the proofs besides those of four proofs of theorems which have many consequences.

In Section \ref{Elementary} we recall a number of elementary definitions of some special subsets of the sphere $S^d$. 
In particular, the notions of hemisphere, lune, $k$-dimensional great subsphere and spherical ball. 
Basic properties of convex bodies on $S^d$ are presented in Section \ref{ConvexBodies}.
We find there the notion of supporting hemisphere of a spherical convex body, as well as the definitions of convex polygons and Reuleaux odd-gons on $S^2$ and rotational convex bodies on $S^d$.
In Section \ref{Lunes} we recall some simple claims on lunes and their corners.
In Section \ref{Width} we review the notions of width and thickness of a spherical convex body $C \subset S^d$ and a number of properties of them.
There is a useful theorem which says how to find the width of $C$ determined by a supporting hemisphere in terms of the radii of the concentric balls touching $C$ from inside or outside.
In Section \ref{Reduced on $S^d$} we find the definition, examples and properties of reduced spherical convex bodies on $S^d$.
First of all a helpful theorem says that through every extreme point $e$ of a reduced body $R \subset S^d$ a lune $L \supset  R$ of thickness $\Delta(R)$ passes with $e$ as the center of one of the two $(d-1)$-dimensional hemispheres bounding $L$. 
There is also a theorem stating that every reduced spherical body of thickness greater than $\frac{\pi}{2}$ is smooth, and a few corollaries.
Section \ref{Reduced on $S^2$} describes the shape of the boundaries of reduced spherical convex bodies on $S^2$ and shows a few consequences.
A theorem and a problem concern the circumradius.
Section \ref{Polygons} is on reduced polygons on $S^2$.
In particular, it appears that every reduced spherical polygon is an odd-gon of thickness at most $\frac{\pi}{2}$. 
A theorem permits to recognize whether a spherical convex polygon is reduced.
Some estimates of the diameter, perimeter and area of reduced polygons are also considered.
Section \ref{Diameter} concerns the diameter of convex bodies, and in particular reduced bodies.  
Especially, we show the relationship between the diameter and the maximum width of a convex body on $S^d$ and we find an estimate for of the diameter of reduced bodies on $S^2$ in terms if their thickness. 
The notion of body of constant diameter is recalled.
Section \ref{Constant Width} is devoted to spherical bodies of constant width. 
We recall that if a reduced body $R \subset S^d$ has thickness at least $\frac{\pi}{2}$, then $R$ is a body of constant width, and
that every body of constant width smaller than $\frac{\pi}{2}$ is strictly convex.
Section \ref{Complete} concerns complete spherical convex bodies.
It appears that complete bodies, constant width bodies and constant diameter bodies on $S^d$ coincide.
Section \ref{Final} compares the notions of spherically convex bodies and their width, as well as the notion of constant width discussed in this chapter, with some earlier analogous notions considered by other authors.

\vskip-0.8cm
{\ }

\section{Elementary notions}\label{Elementary}

\vskip0.11cm
\noindent
Let $S^d$ be the unit sphere in the $(d+1)$-dimensional Euclidean space $E^{d+1}$, where $d\geq 2$. 
The intersection of $S^d$ with an $(m+1)$-dimensional subspace of $E^{d+1}$, where $0 \leq m < d$, is called an {\bf $m$-dimensional great subsphere of $S^d$}\index{great subsphere}. 
In particular, if $m=0$, we get the $0$-dimensional subsphere being {\bf a pair of antipodal points}\index{antipodal point} or {\bf antipodes} in short, and if $m=1$ we obtain the so-called {\bf great circle}\index{great circle}. 
Important for us is the case when $m = d-1$.
In particular, for $S^2$ the $(d-1)$-dimensional great spheres are great circles.

Observe that if two different points are not antipodes, there is exactly one great circle containing them.
If different points $a, b \in S^d$ are not antipodes, by the {\bf arc}\index{arc} $ab$ connecting them we mean the shorter part of the great circle containing $a$ and $b$. 

By the {\bf spherical distance}\index{spherical distance} $|ab|$, or, in short, {\bf distance}\index{distance}, of these points we understand the length of the arc connecting them. 
Moreover, we put $\pi$, if the points are antipodes and $0$ if the points coincide. 
Clearly $|ab| = \angle aob$, where $o$ is the center of $E^{d+1}$.
The {\bf distance from a point $p \in S^d$ to a set}\index{distance to a set} $A \subset S^d$ is understood as the infimum of distances from $p$ to points of $A$.

Consider a non-empty set $A \subset S^d$.
Its interior with respect to the smallest great subsphere of $S^d$ which contains $A$ is called the {\bf relative interior}\index{relative interior}.

The {\bf diameter}\index{diameter} $\diam (A)$ of a set $A \subset S^d$ is the supremum of the spherical distances between pairs of points of $A$.
Moreover, we agree that the empty set and any one-point set have diameter $0$. 
Clearly, if $A$ is closed and has a positive diameter, then the diameter of $A$ is realized for at least one pair of points of $C$.

By a {\bf spherical ball of radius $\rho \in (0, {\pi \over 2}]$}\index{spherical ball}, or shorter, {\bf a ball}\index{ball}, we understand the set of points of $S^d$ having distances at most $\rho$ from a fixed point, called the {\bf center} of this ball. 
An {\bf open spherical ball}\index{open ball}, or shortly {\bf open ball} is the set of points of $S^d$ having distance smaller than $\rho$ from a point.
Balls on $S^2$ are called {\bf disks}\index{disk}. 

Spherical balls of radius $\pi \over 2$ are called {\bf hemispheres}\index{hemisphere}.
In other words, by a {\bf hemisphere} of $S^d$ we mean the common part of $S^d$ with any closed half-space of $E^{d+1}$.
We denote by $H(p)$ the hemisphere whose center is $p$.
Two hemispheres whose centers are antipodes are called {\bf opposite hemispheres}\index{opposite hemispheres}.
By an {\bf open hemisphere}\index{open hemisphere} we mean the set of points having distance less than $\pi \over 2$ from a fixed point.
Hemispheres of $S^d$ play the role of half-spaces of $E^d$. 

By a {\bf spherical $(d-1)$-dimensional ball of radius $\rho \in (0, {\pi \over 2}]$}\index{$(d-1)$-dimensional ball}  we mean the set of points of a $(d-1)$-dimensional great sphere of $S^d$ which are at distance at most $\rho$ from a fixed point, called the {\bf center} of this ball.  
The $(d-1)$-dimensional balls of radius $\pi \over 2$ are called {\bf $(d-1)$-dimensional hemispheres}\index{$(d-1)$-dimensional hemisphere}.
If $d=2$, we call them {\bf semi-circles}\index{semi-circle}.

If an arc $ab$ is a subset of a hemisphere $H$ with $a$ in the boundary $\bd (H)$ of $H$ and $ab$ orthogonal to $\bd (H)$, we say that $ab$ and $\bd (H)$ are orthogonal at $a$, or in short, that they are {\bf orthogonal}\index{orthogonal}.  

We say that a set $C \subset S^d$ is {\bf convex}\index{convex set} if it does not contain any pair of antipodes and if together with every two points of $C$ the whole arc connecting them is a subset of $C$. 

Clearly, for any non-empty convex subset $A$ of $S^d$ there exists a smallest, in the sense of inclusion, non-empty subsphere containing $A$. 
It is unique.

Of course, the intersection of every family of convex sets is also convex. 
Thus for every set $A \subset S^d$ contained in an open hemisphere of $S^d$ there exists a unique smallest convex set containing $A$. 
It is called {\bf the convex hull of}\index {convex hull} $A$ and it is denoted by ${\rm conv} (A)$.
The following fact from \cite{L3} results by applying an analogous theorem for compact sets in $E^{d+1}$.

\begin{cla} \label{conv-closed} 
If $A \subset S^d$ is a closed subset of an open hemisphere, then ${\rm conv} (A)$ is also closed. 
\end{cla}

By the well-known fact that a set $C \subset S^d$ is convex if and only if the cone generated by it in $E^{d+1}$ is convex, and by the separation theorem for convex cones in $E^{d+1}$ we obtain, as observed in \cite{L3}, the following analogous fact for $S^d$.
By the way, the statement is later derived also by Z\u alinescu \cite{Za} from his Theorem 2. 

\begin{cla} \label{opposite-hemispheres} 
Every two convex sets on $S^d$ with disjoint interiors are subsets of two opposite hemispheres.
\end{cla}

\vskip0.21cm

\section{Convex bodies} \label{ConvexBodies}

\vskip0.11cm
By a {\bf convex body}\index{convex body} on $S^d$ we mean a closed convex set with non-empty interior. 

If a $(d-1)$-dimensional great sphere $G$ of $S^d$ has a common point $t$ with a convex body $C \subset S^d$ and if its intersection with the interior of $C$ is empty, we say that $G$ is a {\bf supporting $(d-1)$-dimensional great sphere of $C$ passing through $t$}\index{supporting great sphere}. 
We also say that $G$ {\bf supports} $C$ at $t$.
If $H$ is the hemisphere bounded by $G$ and containing $C$, we say that $H$ is a {\bf supporting hemisphere of $C$}\index{supporting hemisphere} and that $H$ {\bf supports $C$ at $t$}.

If at a boundary point $p$ of a convex body $C \subset S^d$ exactly one hemisphere supports $C$, we say that $p$ is a {\bf smooth point}\index{smooth point} of the boundary $\bd (C)$ of $C$.
If all boundary points of $C$ are smooth, then $C$ is called a {\bf smooth body}\index{sooth body}.

Let $C \subset S^d$ be a convex body and let $Q \subset S^d$ be a convex body or a hemisphere.   
We say that $C$ {\bf touches $Q$ from outside}\index{touch from outside} if $C \cap Q \not = \emptyset$ and ${\rm int} (C) \cap {\rm int} (Q) = \emptyset$. 
We say that $C$ {\bf touches $Q$ from inside}\index{touch from inside} if $C \subset Q$ and $\bd (C) \cap \bd (Q) \not = \emptyset$. 
In both cases, points of $\bd (C) \cap \bd (Q)$ are called {\bf points of touching}\index{point of touching}.

The convex hull $V$ of $k \geq 3$ points on $S^2$ such that each of them does not belong to the convex hull of the remaining points is named a {\bf spherical convex $k$-gon}\index{spherical polygon}. 
The points mentioned here are called the {\bf vertices}\index{vertex} of $V$. 
We write $V= v_1v_2\dots v_k$ provided $v_1, v_2, \dots , v_k$ are successive vertices of $V$ when we go around $V$ on the boundary of $V$. 
In particular, when we take $k\geq 3$ successive points in a spherical circle of radius less than $\pi \over 2$ on $S^2$ with equal distances of every two successive points, we obtain a {\bf regular spherical $k$-gon}\index{regular spherical $k$-gon}. 

Take a regular spherical $k$-gon $v_1v_2\dots v_k \subset S^2$, where $k \geq 3$ is odd.
Clearly, all the distances $|v_iv_{i+ {{k-1} \over 2}}|$ and $|v_iv_{i+ {{k+1} \over 2}}|$ for $i=1, \dots ,k$ are equal (the indices are taken modulo $k$). 
Denote them by $\delta$.
Assume that $\delta \leq {\pi \over 2}$. 
Let $B_i$, where $i=1,\dots ,k$, be the disk with center $v_i$ and radius $\delta$.
The set $B_1 \cap \dots \cap B_k$ is called a {\bf spherical Reuleaux $k$-gon}\index{Reuleaux polygon}.
Clearly, it is a convex body.

We say that $e$ is an {\bf extreme point}\index{extreme point} of a convex body $C \subset S^d$ provided the set $C \setminus \{ e \}$ is convex.
From the analogue of the Krein-Milman theorem for convex cones (e.g., see \cite{FL}) its analogue for spherical convex bodies formulated on p. 565 of \cite{L3} follows.

\begin{cla}
Every convex body $C \subset S^d$ is the convex hull of its extreme points.
\end{cla}

This and the fact that the intersection of any closed convex body $C \subset S^d$ with any of its supporting $(d-1)$-dimensional great spheres is a closed convex set imply the following fact from \cite{L3}.

\begin{cla} \label{passes-extreme}  
The boundary of every supporting hemisphere of a convex body $C \subset S^d$ passes through an extreme point of $C$.
\end{cla}

The next property is shown in \cite{LaMu2}.

\begin{cla}\label{extreme} 
Let $C\subset S^d$ be a convex body. 
Every point of $C$ belongs to the convex hull of at most $d+1$ extreme points of $C$.
\end{cla}

Let us add that Shao and Guo \cite{SG} proved analytically that every closed set $C\subset S^d$ is the convex hull of its extreme points.

We say that two sets on $S^2$ are {\bf symmetric with respect to a great circle}\index{symmetric} if they are symmetric with respect to the plane of $E^3$ containing this circle.
More generally, two sets on $S^d$ are called {\bf symmetric with respect to a $(d-1)$-dimensional subsphere $S^{d-1}$} if they are symmetric with respect to the hyperplane of $E^{d+1}$ containing $S^{d-1}$.
If a set coincides with its symmetric, we say that the set is {\bf symmetric with respect to a great circle}.

Let us construct a rotational body on $S^3$. 
Take a two-dimensional subsphere $S^2$ of $S^3$ and a convex body $C$ symmetric with respect to a great circle of $S^2$. 
Denote by $A$ the common part of this great circle and $C$, and call it an {\bf arc of symmetry of $C$}\index{arc of symmetry}. 
Take all subspheres $S_{\lambda}^2$ of $S^3$ containing $A$, and on each $S_{\lambda}^2$ take the copy $C_{\lambda}$ of $C$ at the same place with respect to $A$. 
The union ${\rm rot}_A (C)$ of all these $C_{\lambda}$ (including $C$) is called a {\bf rotational body}\index{rotational body} on $S^3$, or the {\bf body of rotation of $C$ around $A$}\index{rotation around an arc}. 
An analogous construction can be provided on $S^{d+1}$ in place of $S^3$; this time we ``rotate"  a convex body $C \subset S^d$ having a $(d-1)$-dimensional subsphere $S^{d-1}$ of symmetry. 
We rotate it around the common part $A$ of $S^{d-1}$ and $C$.
Again we denote the obtained body by ${\rm rot}_A (C)$.

Since $S^d$ with the distance of points defined in the preceding section is a metric space, we may consider the Hausdorff distance\index{Hausdorff distance} on $S^d$. 
Applying the classical Blaschke selection theorem\index{Blaschke selection theorem} (see e.g. p. 64 of the monograph \cite{Eg} by Eggleston and p. 50 of the monograph by Schneider \cite{Sch}) in $E^{d+1}$ we easily obtain the following spherical version of Blaschke spherical theorem mentioned in \cite{L1} before Lemma 4. 
See also the paper of Hai and An \cite{HA}, where Corollary 4.11 says that {\it if $X$ is a proper uniquely geodesic space, then  every uniformly bounded sequence of nonempty geodesically convex subsets of $X$ contains a subsequence which converges to a nonempty compact geodesically convex subset of $X$}. 
For the context see the book \cite{Pap1} by Papadopoulos. 
Moreover look to Theorem 3.2 of the paper \cite{BoS} by B\"oroczky and Sagemeister.

\begin{thm}
Every sequence of convex bodies on $S^d$ contains a subsequence convergent to a convex body on $S^d$.
\end{thm}

\begin{cor} \label{sequence1} 
Every sequence of convex bodies being subsets of a fixed convex body $C \subset S^d$ contains a subsequence convergent to a convex body being a subset of $C$.
\end{cor}

\begin{cor} \label{circumscribed} 
For every convex body $C \subset S^d$ there exist a ball of minimum radius containing $C$ and a ball of maximum radius contained in $C$.
\end{cor}

Proofs of these corrolaries are similar to proofs of analogous facts in $E^d$. 
For instance, see 
Section 18 of the book \cite{Ben} by Benson and Section 16 of the book \cite{Lay} by Lay.

\vskip0.11cm

\section{Lunes}\label{Lunes}

\vskip0.11cm
\noindent
Recall the classical notion of lune.
For any distinct non-opposite hemispheres $G$ and $H$ of $S^d$ the set $L = G \cap H$ is called {\bf a lune}\index{lune} of $S^d$. 
Lunes of $S^d$ play the role of strips in $E^d$.
The $(d-1)$-dimensional hemispheres bounding the lune $L$ which are contained in $G$ and $H$, respectively, are denoted by $G/H$  and $H/G$.
Here is a claim from \cite{L3}.

\begin{cla} \label{pair-determines} 
Every pair of different points  $a,b \in S^d$ which are not antipodes determines exactly one lune $L$ such that  $a$,  $b$ are the centers of the $(d-1)$-dimensional hemispheres bounding~$L$. 
\end{cla}

Since every lune $L$ determines exactly one pair of centers of the $(d-1)$-dimensional hemispheres bounding $L$, from Claim \ref{pair-determines} we see that there is a one-to-one correspondence between lunes and pairs of different points, which are not antipodes, of $S^d$. 

Observe that $(G/H) \cup (H/G)$ is the boundary of the lune $G \cap H$ (in particular, every lune of $S^2$ is bounded by two different semi-circles).  
Denote by $c_{G/H}$ and $c_{H/G}$ the centers of $G/H$ and $H/G$, respectively.
By {\bf corners}\index{corner} of the lune $G \cap H$ we mean points of the set $(G/H) \cap (H/G)$. 
Of course, $r \in (G/H) \cup (H/G)$ is a corner of $G \cap H$ if and only if $r$ is equidistant from $c_{G/H}$ and $c_{H/G}$.
Observe that the set of corners of a lune is a $(d-2)$-dimensional subsphere of $S^d$.
In particular, every lune of $S^2$ has two corners. 
They are antipodes.

By the {\bf thickness $\Delta (L)$ of a lune}\index{thickness of lune} $L = G \cap H \subset S^d$ we mean the spherical distance of the centers of the $(d-1)$-dimensional hemispheres $G/H$ and $H/G$ bounding $L$.
Observe that it is equal to each of the non-oriented angles 
$\angle c_{G/H}rc_{H/G}$, where $r$ is any corner of $L$. 

Here is a fact shown in \cite{L3}.

\begin{cla} \label{DistInLunes} 
Let $H$ and $G$ be different and not opposite hemispheres. 
Consider the lune $L= H \cap G$. 
Let $x \not = c_{G/H}$ belong to $G/H$. 
If $\Delta (L) < {\pi \over 2}$, we have $|xc_{H/G}| > |c_{G/H}c_{H/G}|$.
If $\Delta (L) = {\pi \over 2}$, we have $|xc_{H/G}| = |c_{G/H}c_{H/G}|$.
If $\Delta (L) > {\pi \over 2}$, we have $|xc_{H/G}| < |c_{G/H}c_{H/G}|$.
\end{cla}

Similarly to Claim \ref{sequence1} we obtain the following Claim formulated in \cite{L3} and needed for the proof of Theorem \ref{e-reduced}.

\begin{cla} \label{sequence2} 
Every sequence of lunes of a fixed thickness $t$ of $S^d$ contains a subsequence of lunes convergent to a lune of thickness $t$.
\end{cla}

The following two claims are from \cite{L5}.
 
\begin{cla}  \label{twolunes} 
Consider a hemisphere $H(c)$ of $S^d$.
Then any $(d-1)$-dimensional subsphere of $S^d$ containing $c$ dissects $H(c)$ into two lunes of thickness $\frac{\pi}{2}$.
\end{cla}

\begin{cla}  \label{corners}
Let $L \subset S^d$ be a lune and let $C \subset L$ be a convex body such that the set $F= C \cap {\rm corn} (L)$ is non-empty. 
Then at least one extreme point of $C$ is a corner point of $L$.
\end{cla}

The next claim and remark are taken from \cite{L3}. 

\begin{cla} \label{C-subset-L} 
Let ${\rm  diam} (C) \leq {\pi \over 2}$ for a convex body $C \subset S^d$ and assume that ${\rm  diam} (C) = |ab|$ for some points $a, b \in C$. 
Denote by $L$ the lune such that $a$ and $b$ are the centers of the $(d-1)$-dimensional hemispheres bounding $L$.
We have $C \subset L$.
\end{cla}

In general, this claim does not hold without the assumption that ${\rm  diam} (C) \leq {\pi \over 2}$. 
A simple counterexample is the triangle $T = abc$ with $|ab|={2\over 3}\pi \approx 2.0944$, $|bc| = {\pi \over 6} \approx 0.5236$ and $\angle abc =95{^\circ}$.
From the law of cosines for sides, called also Al'Battani formulas in \cite{Mu}, p. 45,, we get $|ac| \approx 2.0609$. 
Consequently, $|ab| ={2\over 3}\pi$ is the diameter of $T$. 
Since $\angle abc =95{^\circ} $, the lune with centers $a$ and $b$ of the semi-circles bounding it does not contain $c$.
Still its thickness is ${2\over 3}\pi$.
Thus this lune does not contain~$T$.

The subsequent claim is shown in \cite{L6}.

\begin{cla} 
Let $K$ be a hemisphere of $S^d$ and let $p \in \bd (K)$.
Moreover, let $pq \subset K$ be an arc orthogonal to $\bd (K)$ with $q$ in the interior of $K$ and $|pq| < \frac{\pi}{2}$. 
Then amongst all the lunes of the form $K \cap M$, with $q$ in the boundary of the hemisphere $M$, only the lune $K \cap K_\dashv$ such that $pq$ is orthogonal to $\bd (K_\dashv)$ at $q$ has the smallest thickness.
\end{cla}

\vskip0.21cm

\section{Width and thickness of a convex body} \label{Width}

\vskip0.11cm
We say that {\bf a lune passes through a boundary point $p$ of a convex body $C \subset S^d$} if the lune contains $C$ and if the boundary of the lune contains $p$. 
If the centers of both the $(d-1)$-dimensional hemispheres bounding a lune belong to $C$, then we call such a lune an {\bf orthogonally supporting lune of}\index{supporting lune}~$C$. 

For any convex body $C \subset S^d$ and any hemisphere $K$ supporting $C$ we define 
{\bf the width of $C$ determined by $K$}\index{width} as the minimum thickness of a lune $K \cap K^*$ over all hemispheres $K^* \not = K$ supporting $C$.
We denote it by ${\rm width}_K (C)$. 
By compactness arguments we immediately see that at least one such hemisphere $K^*$ exists, and thus at least one corresponding lune $K \cap K^*$ exists.

This notion of width of $C \subset S^d$ is an analogue of the notion of width of a convex body $C \subset E^d$ between a pair of parallel supporting hyperplanes $H_1, H_2$ supporting $C$.
Our $K \subset S^d$ takes over the role of one of the half-spaces bounding $H_1$ or $H_2$ and containing $C$ and our $K^*$ takes over the role of the other half-space, which in $E^d$ is unique.

The following lemma from \cite{L3} is needed for the proof of the forthcoming Theorem~\ref{I-III}.

\vskip0.2cm
\begin{figure}[htbp]        

\includegraphics[width=2.51in]{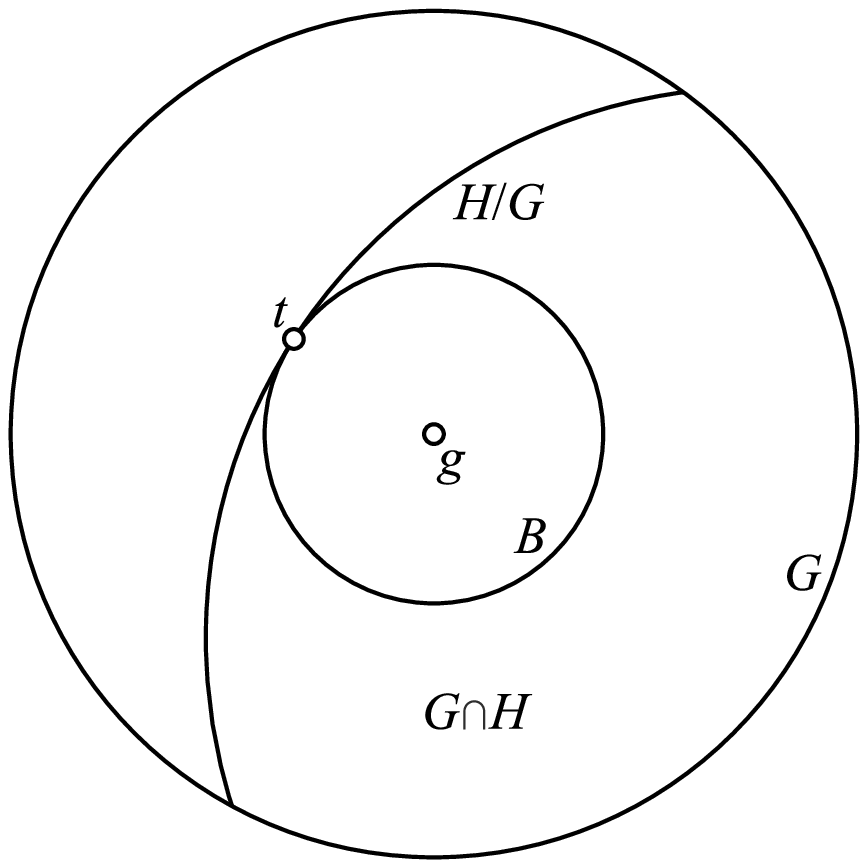} \hskip0.1cm
\includegraphics[width=2.51in]{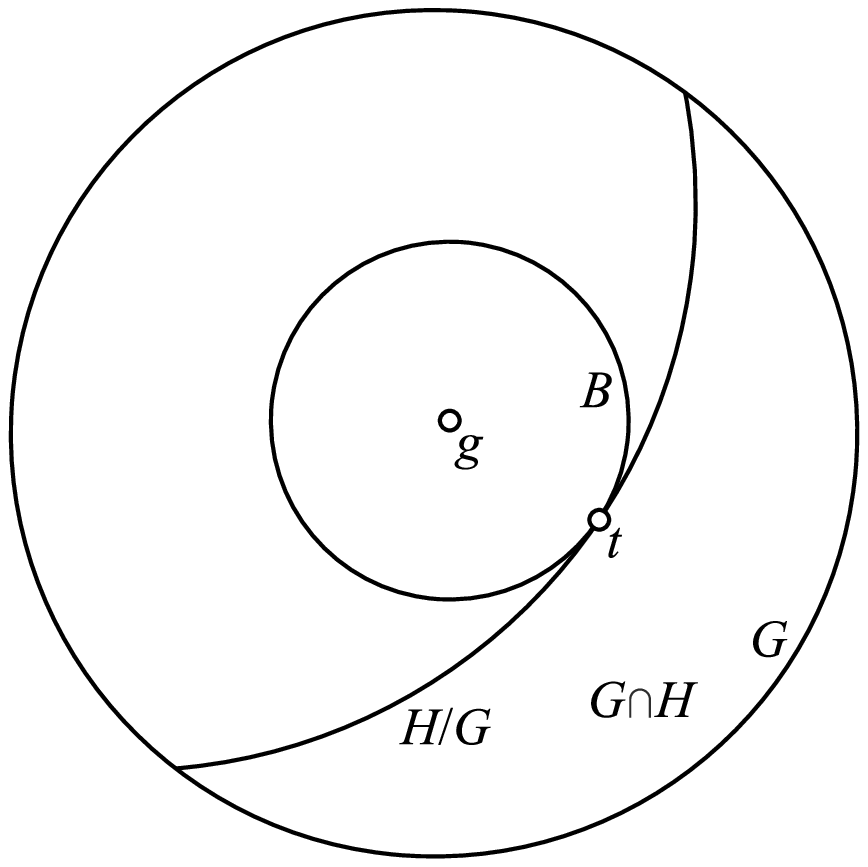}  \\  

\vskip 0.1cm
\caption{Illustration to Lemma \ref{t-is-center-of-H/G}}   \label{fig:touches}

\end{figure}
\vskip0.2cm

\begin{lem} \label{t-is-center-of-H/G} 
Let $G$ and $H$ be different and not opposite hemispheres, and let $g$ denote the center of $G$. 
If $g \not \in \bd (H)$, then by $B$ denote the ball with center $g$ which touches $H$ (from inside or outside) and by $t$ the point of touching.
If $g \in \bd (H)$, we put $t=g$. 
We claim that $t$ is always at the center of the $(d-1)$-dimensional hemisphere $H/G$. 
\end{lem}

Figure \ref{fig:touches} shows the two situations when $B$ touches $H$ from inside and outside as in the text of the lemma.
Here and in later figures we adopt an orthogonal look to the sphere from outside. 

The theorem from \cite{L3} given below and its proof present a procedure for establishing 
the width ${\rm width}_K (C)$ of a convex body $C \subset S^d$ in terms of the radii of balls concentric with $K$ and touching $C$.

\begin{thm} \label{I-III} 
Let $K$ be a hemisphere supporting a convex body $C \subset S^d$. 
Denote by $k$ the center of $K$.

\vskip0.1cm
I. If $k \not \in C$, then there exists a unique hemisphere $K^*$ supporting $C$ such that the lune $L= K \cap K^*$ contains $C$ and has thickness ${\rm width}_K (C)$.  
This hemisphere supports $C$ at the point $t$ at which the largest ball $B$ with center $k$ touches $C$ from outside.
We have $\Delta (K \cap K^*) = {\pi \over 2} - \rho_B$, where $\rho_B$ denotes the radius of $B$.

\vskip0.1cm
II. If $k \in \bd (C)$, then there exists at least one hemisphere $K^*$ supporting $C$ such that $L= K \cap K^*$ is a lune containing $C$ which has thickness  ${\rm width}_K (C)$. 
This hemisphere supports $C$ at $t=k$. 
We have $\Delta (K \cap K^*) = {\pi \over 2}$.

\vskip0.1cm
III. If $k \in {\rm int} (C)$, then there exists at least one hemisphere $K^*$ supporting $C$ such that $L= K \cap K^*$ is a lune containing $C$ which has thickness  ${\rm width}_K (C)$.  
Every such $K^*$ supports $C$ at exactly one point $t \in \bd (C) \cap B$, where $B$ denotes the largest ball with center $k$ contained in $C$, and for every such $t$ this hemisphere $K^*$, denoted by $K_t^*$, is unique.
For every $t$ we have $\Delta (K \cap K_t^*) = {\pi \over 2} + \rho_B$, where $\rho_B$ denotes the radius of $B$.  
\end{thm}

\begin{proof}
In Figure \ref{fig:I-III} we see two illustrations to this theorem and its proof.  
They, respectively, concern Parts I and III.

\vskip0.1cm
\begin{figure}[htbp]        

\includegraphics[width=2.51in]{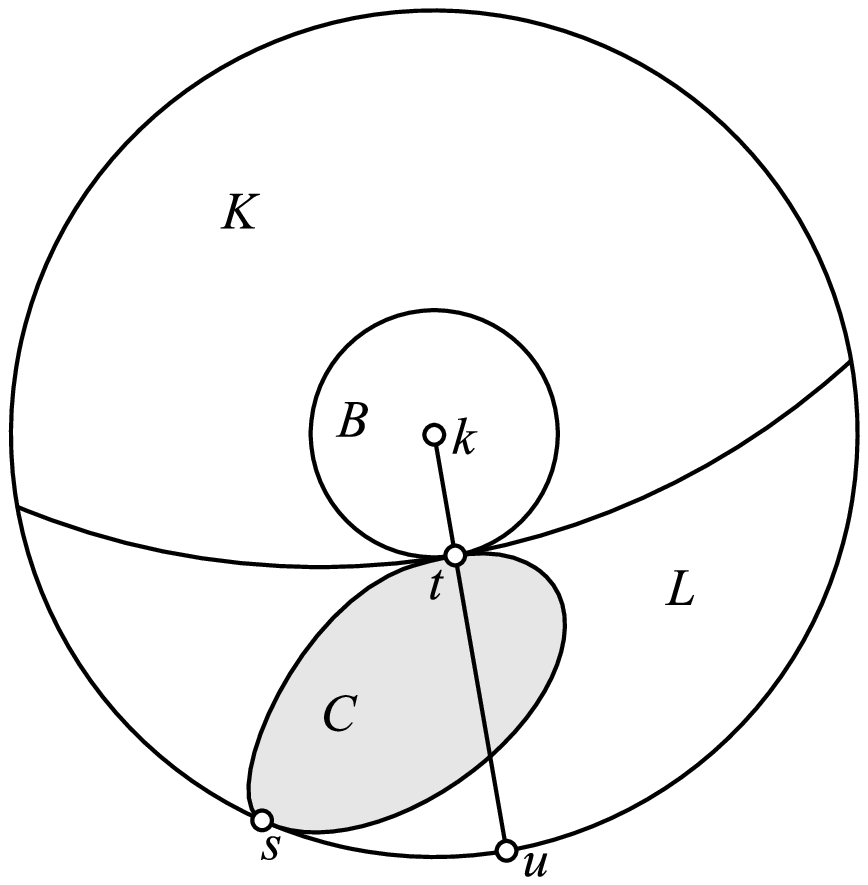} \hskip0.1cm 
\includegraphics[width=2.51in]{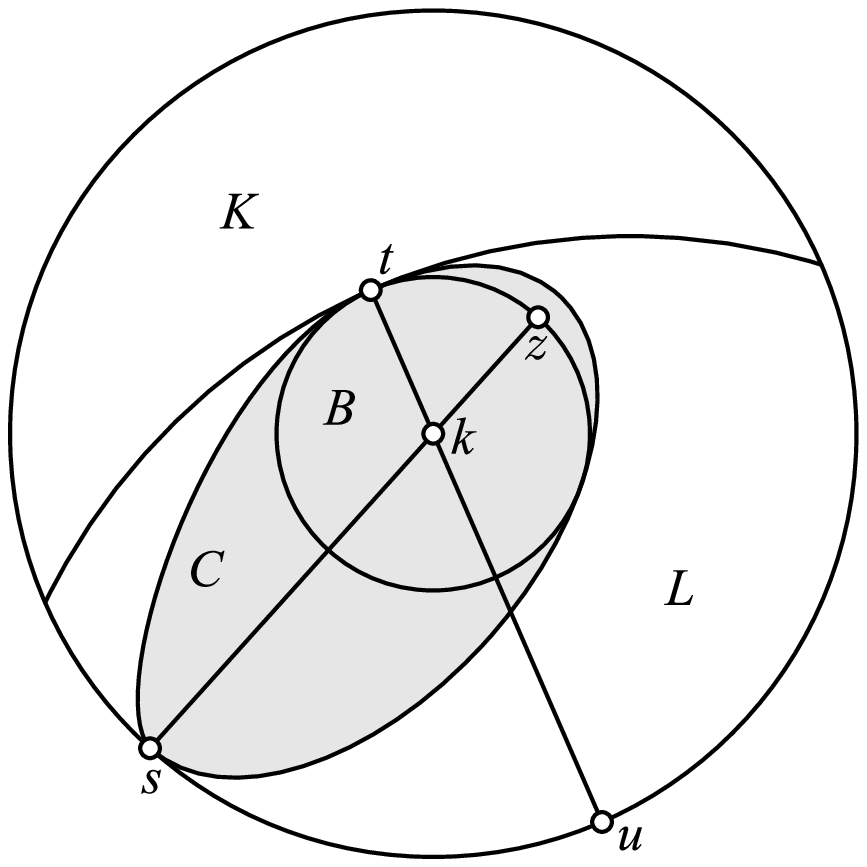} \\ 

\caption{Illustration to Theorem 1 and its proof}   \label{fig:I-III}

\end{figure}

\vskip0.1cm
Part I. 
Since $C$ is a convex body and $B$ is a ball, we see that $B$ touches $C$ from the outside and the point of touching is unique.
Denote it by $t$ (see the first illustration of Figure \ref{fig:I-III}). 
By Lemma \ref{opposite-hemispheres}, the bodies $C$ and $B$ are in some two opposite hemispheres.
What is more, since $B$ is a ball touching $C$ from the outside, 
this pair of hemispheres is unique. 
Denote by $K^*_t$ the one of them which contains $C$.
We intend to show that $K^*_t$ is nothing else but the promised $K^*$.

Denote by $k^*$ the center of $K^*_t$.
Since $k$ is also the center of $B$ and since $B$ and $K^*_t$ touch from the outside at $t$, we have $t \in kk^*$. 
From Lemma 5 we see that $t$ is the center of the $(d-1)$-dimensional hemisphere $K^*_t/K$. 
Analogously, from this lemma we conclude that the common point $u$ of $kk^*$ and the boundary of $K$ is the center of $K/{K^*_t}$. 
Since $t$ and $u$ are centers of the $(d-1)$-dimensional hemispheres bounding the lune $K \cap K^*_t$, we have $|tu| = \Delta (K \cap K^*_t)$.
This and $|kt| + |tu| = |ku| = {\pi \over 2}$ imply $\Delta (K \cap K^*_t) = {\pi \over 2} - \rho_B$.

If we assume that there exists a hemisphere $M \supset C$ with $\Delta (K \cap M) < {\pi \over 2} - \rho_B$, then the lune $K \cap M$ must be disjoint from $B$, and hence it does not contain $C$.
A contradiction.
Thus $K \cap K^*_t$ is a narrowest lune of the form $K \cap N$ containing $C$.
It is the unique lune of this form in virtue of the uniqueness of $t$ and $K^*_t$ explained at the beginning of the proof of Part I

\vskip0.1cm
Part II. 
Clearly, there is at least one hemisphere $K^*$ supporting $C$ at $k$. 
Of course, $\Delta (K \cap K^*) = {\pi \over 2}$.
By Lemma \ref{t-is-center-of-H/G} we see that $k$ is the center of $K^*/K$.  

\vskip0.1cm
Part III. 
Take the largest ball $B \subset C$ with center $k$. 
Clearly, there is at least one boundary point $t$ of $C$ which is also a boundary point of $B$ (see the second illustration of Figure \ref{fig:I-III}). 
We find a hemisphere $K^*_t$ which supports $C$ at $t$.
Of course, it also supports $B$ and thus, for given $t$, it is unique.  

\vskip0.1cm
For every $t$ there is a unique point $u \in K/{K^*_t}$ such that $k \in tu$.
Thus, $|ku| = {\pi \over 2}$ and $|kt| = \rho_B$ imply $|tu| = {\pi \over 2} + \rho_B$. 
Hence the facts, resulting from Lemma \ref{t-is-center-of-H/G}, that $t$ is the center of $K^*_t/K$ and that $u$ is the center of $K/{K^*_t}$ give $\Delta (K \cap K^*_t) = {\pi \over 2} + \rho_B$. 

If we assume that there exists a hemisphere $M \supset C$ such that the lune $K \cap M$ is narrower than ${\pi \over 2} + \rho_B$, then this lune does not contain $B$, and hence it also does not contain $C$, a contradiction.
Thus the narrowest lunes of the form $K \cap N$ containing $C$ are of the form $K \cap K^*_t$.
\end{proof}

Let us point out that in Part I, if the center $k$ of $K$ does not belong to $C$, the lune $K \cap K^*$ is unique. 
In Part II this narrowest lune $K \cap K^*$ containing $C$ sometimes is unique and sometimes not. 
This depends on the point $k=t$ of $C$ which belongs to the boundary of $B$.
In Part III for any given point $t$ of touching $C$ by $B$ from the inside (we may have one, or finitely many, or infinitely many such points $t$), the lune $K \cap K^*_t$ is unique.

For instance, if $C \subset S^2$ is a regular spherical triangle of sides $\pi \over 2$ and the circle bounding a hemisphere $K$ contains a side of this triangle, then $K \cap K^*$ is not unique.
Namely, as $K^*$ we may take any hemisphere containing $C$, whose boundary contains this vertex of $C$ which does not belong to $K$. 
The thickness of every such lune $K \cap K^*$ is equal to $\pi \over 2$.
If $C$ is a regular spherical triangle of sides greater than $\pi \over 2$ and the boundary of $K$ contains a side of this triangle, then $K \cap K^*$ is also not unique.
This time the boundary of $K^*$ contains a side of $C$ different from the side which is contained in $K$. 
So we have exactly two positions of $K^*$.

Here are two corollaries from Theorem \ref{I-III}. 
For the second we additionally apply Lemma \ref{t-is-center-of-H/G}.

\begin{cor} \label{Cor1LasWIDTH} 
If $k \not \in C$, then ${\rm width}_K (C) = {\pi \over 2} - \rho_B$.
If $k \in \bd (C)$, we have ${\rm width}_K (C) = {\pi \over 2}$.
If $k \in {\rm int} (C)$, then ${\rm width}_K (C) = {\pi \over 2} + \rho_B$.
\end{cor}

\begin{cor} \label{Cor2LasWIDTH} 
The point $t$ of support in Theorem \ref{I-III} is the center of the $(d-1)$-dimensional hemisphere $K^*/K$. 
\end{cor}

The following theorem is also proved in \cite{L3}.

\begin{thm} 
As the position of the $(d-1)$-dimensional supporting hemisphere of a convex body $C \subset S^d$ changes, the width of $C$ determined by this hemisphere changes continuously. 
\end{thm}

\noindent
\begin{proof} 
We keep the notation of Theorem \ref{I-III}.
Of course, the positions of $k$ and thus that of $B$ depend continuously on $K$.
Hence ${\pi \over 2} - \rho_B$ and ${\pi \over 2} + \rho_B$ change continuously.
This fact and Corollary \ref{Cor1LasWIDTH} imply the thesis of our theorem.
It does not matter here that for a fixed $K$ sometimes the lunes $K \cap K^*$ are not unique; still they are all of equal thickness.
\end{proof}

Compactness arguments lead to the conclusion that for every convex body $C \subset S^d$ the supremum of ${\rm width}_H (C)$ over all hemispheres $H$ supporting $C$ is realized for a supporting hemisphere of $C$, that is, we may take here the maximum instead of supremum.

We define the {\bf thickness $\Delta (C)$ of a convex body}\index{thickness} $C \subset S^d$ as the minimum of the thickness of the lunes containing $C$ (by compactness arguments, this minimum is realized), or even more intuitively as the thickness of each ``narrowest'' lune containing $C$.
In other words,

\vskip-0.5cm
$$\Delta (C) =\min \{ {\rm width}_K (C); K {\rm \ is \ a \ supporting \ hemisphere \ of} \ C \}.$$

\begin{exa} \label{ExRegTria}
{\rm Applying Theorem \ref{I-III} we easily find the thickness of any regular triangle $T_\alpha \subset S^2$ of angles $\alpha$.
Formulas of spherical trigonometry imply that $\Delta (T_\alpha) = {\rm arc cos} {\cos \alpha \over \sin \alpha /2}$ for $\alpha < {\pi \over 2}$. 
If $\alpha \geq {\pi \over 2}$ (but, of course, $\alpha < {2\over 3}\pi$), then $\Delta (T_\alpha) = \alpha$. 
In both cases $\Delta (T_\alpha)$ is realized for ${\rm width}_K (T_\alpha)$, where $K$ is a hemisphere whose bounding semicircle contains a side $S$ of $T_\alpha$.
In the first case $T_\alpha$ is symmetric with respect to the great circle containing the arc $kk^*$ connecting the centers $k$ and $k^*$ of $K/K^*$ and $K^*/K$, respectively, while in the second  $T_\alpha$ is symmetric with respect to the great circle passing through the middle of $kk^*$ and having endpoints at the corners of the lune $K \cap K^*$ (observe that this time we have exactly two positions of $K^*$, each with $\bd (K^*)$ containing a side of $T_\alpha$ different from $S$). 
For $\alpha = {\pi \over 2}$ there are infinitely many positions in $K^*$; just these hemispheres which support $T_\alpha$ at $k^*$.} 
\end{exa}

\vskip0.2cm
The first thesis of the following fact was first proved in  \cite{L3}. 
The last statement was given in \cite{LaMu1} for $S^2$ and in \cite{LaMu2} for $S^d$. 

\begin{cla} \label{BothTheCenters} 
Consider a convex body $C \subset S^d $ and any lune $L$ of thickness $\Delta (C)$ containing~$C$. 
Both centers of the $(d-1)$-dimensional hemispheres bounding $L$ belong to $C$.
Moreover, if $\Delta (C) > \frac{\pi}{2}$, then both centers are smooth points of the boundary of $C$.
\end{cla}

\vskip0.21cm

\section{Reduced spherical bodies} \label{Reduced on $S^d$}

\vskip0.11cm
\noindent
Reduced bodies in the Euclidean space $E^d$ are defined by Heil \cite{He} and next considered in very many papers (references of many of them can be found in the survey articles \cite{LaMa1} and \cite{LaMa2}).
In analogy to this notion we define reduced convex bodies on $S^d$. 
Namely, after \cite{L3} we say that a convex body $R \subset S^d$ is {\bf reduced}\index{reduced body} if $\Delta (Z) < \Delta (R)$ for every convex body $Z \subset R$ different from $R$. 

It is easy to show that all regular odd-gons on $S^2$ of thickness at most $\pi \over 2$ are reduced bodies. 
The assumption that the thickness is at most  $\pi \over 2$ matters here.
For instance take the regular triangle $T_\alpha$ of angles $\alpha > {\pi \over 2}$ (see Example \ref{ExRegTria}). 
Take the hemisphere $K$ whose boundary contains a side of $T_\alpha$ and apply Part III of Theorem \ref{I-III}. 
The corresponding ball $B \subset T_\alpha$ touches $T_\alpha$ from the inside at exactly two points $t_1, t_2$. 
Cutting off a part of $T_\alpha$ by the shorter arc of the boundary of $B$ between $t_1$  and $t_2$ we obtain a convex body $Z \subset T_\alpha$. 
We have $\Delta (Z) = \Delta (T_\alpha)$, which implies that $T_\alpha$ is not reduced.

Of course, the spherical bodies of constant width on $S^d$ considered later in Section \ref{Constant Width} are reduced bodies.
In particular, every Reuleaux polygon is a reduced body on $S^2$. 

Dissect a disk on $S^2$ by two orthogonal great circles through its center. 
The four parts obtained are called {\bf quarters of a disk}\index{quarter of disk}. 
In particular, the triangle of sides and angles $\pi \over 2$ is a quarter of a disk.
It is easy to see that every quarter of a disk is a reduced body and that its thickness is equal to the radius of the original disk.
More generally, each of the $2^d$ parts of a spherical ball on $S^d$ dissected by $d$ pairwise orthogonal great $(d-1)$-dimensional spheres through the center of this ball is a reduced body on $S^d$. 
We call it a ${1 \over 2^d}${\bf -th part of a ball}\index{part of a ball}.
Clearly, its thickness is equal to the radius of the above ball.

\vskip0.15cm
\begin{figure}[htbp]
\hskip 3cm \includegraphics[width=6.77cm,height=7.55cm]{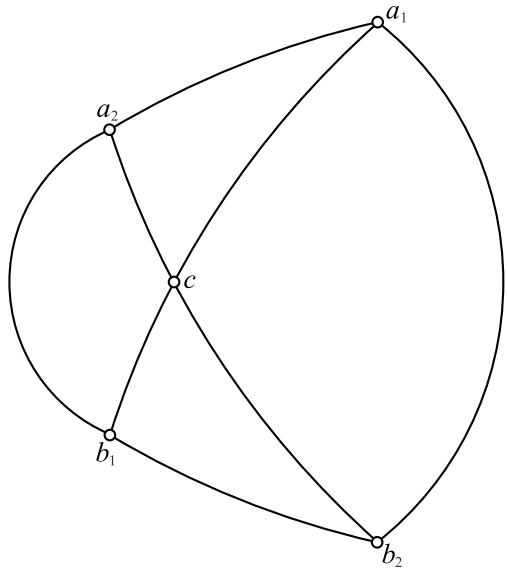} \\ 
\vskip-0.2cm
\caption{The reduced body on $S^2$ from Example \ref{ExRedBodyOn$S^2$}} \label{fig:ExReduced}
\end{figure}

\vskip-0.1cm
\begin{exa} \label{ExRedBodyOn$S^2$} 
{\rm On $S^2$ consider arcs $a_1b_1$ and $a_2b_2$ of equal lengths and with a common point $c$ such that $ca_2 \perp a_2a_1$,  $cb_1 \perp b_1b_2$ and $|ca_2| = |cb_1|$.
As a consequence,  $|ca_1| \geq |cb_1|$ and $|cb_1| \leq |cb_2|$. 
Provide the shorter piece $P_1$ of the circle with center $c$ connecting $a_1$ and $b_2$. 
Also provide the shorter piece $P_2$ of the circle with center $c$ connecting $a_2$ and $b_1$. 
Let $P = \conv(P_1 \cup P_2)$ (see Figure \ref{fig:ExReduced}).
Observe that $P$ is a reduced body of thickness $|a_1b_1| = |a_2b_2|$.}
\end{exa}

\vskip0.2cm
Consider a reduced body $R \subset S^{d-1}$ having a $(d-2)$-dimensional subsphere $S^{d-2}$ of symmetry. 
Put $A = S^{d-2} \cap R$.
Then ${\rm rot}_A(R)$ is a reduced body on $S^d$.
For instance, when we take as $R \subset S^2$ a triangle with an axis of symmetry of it, we get a ``spherical cone" on $S^3$, which is a spherical reduced body.
We can also rotate the reduced body $P$ from Example \ref{ExRedBodyOn$S^2$} taking as $A$ its axis of symmetry, again obtaining a reduced body on $S^3$

The following two theorems are proved in  \cite{L3}.
Having in mind the intuitive ideas of the first proof, we recall it here.

\begin{thm} \label{e-reduced} 
Through every extreme point $e$ of a reduced body $R \subset S^d$ a lune $L \supset  R$ of thickness $\Delta(R)$ passes with $e$ as the center of one of the two $(d-1)$-dimensional hemispheres bounding $L$. 
\end{thm}

\noindent
\begin{proof}
Let $B_i$ be the open ball of radius $\Delta (R) /i$ centered at $e$ and let $R_i = {\rm conv} (R \setminus  B_i)$ for $i=2,3, \dots$ (see Figure \ref{fig:extreme}).
By Claim \ref{conv-closed} every $R_i$ is a convex body.
Moreover, since $e$ is an extreme point of $R$, then $R_i$ is a proper subset of $R$. 
So, since $R$ is reduced, $\Delta(R_i) < \Delta (R)$. 
By the definition of thickness of a convex body, $R_i$ is contained in a lune $L_i$ of thickness $\Delta(R_i)$.

\vskip0.2cm
\begin{figure}[htbp]
\hskip1cm \includegraphics[width=10.83cm,height=8.455cm]{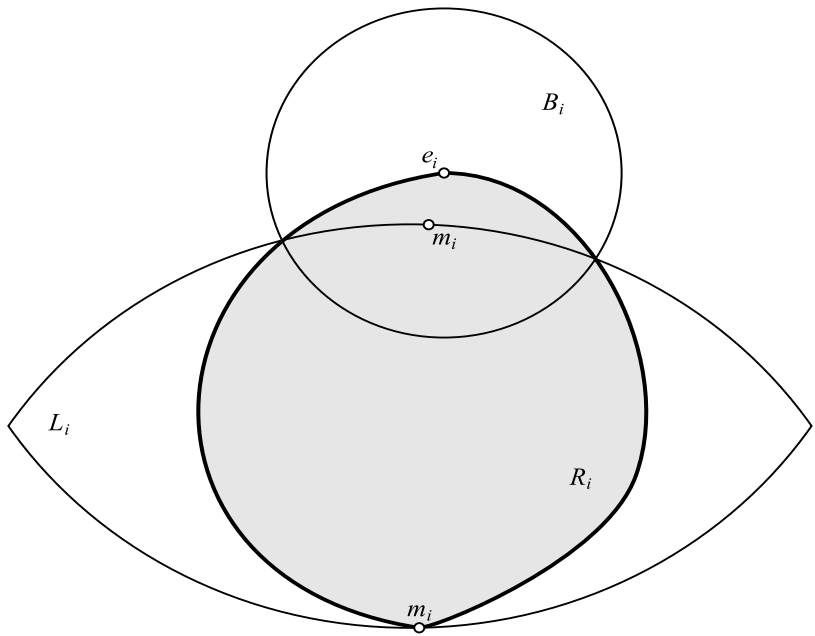} \\ 
\caption{Illustration to the proof of Theorem \ref{e-reduced}} \label{fig:extreme}
\end{figure}

From Claim \ref{sequence2} we conclude that there exists a subsequence of the sequence $L_2, L_3, \dots$ converging to a lune $L$.
Since $R_i \subset L_i$  for $i=2,3, \dots$, we obtain that $R \subset L$. 
Since $\Delta(L_i) = \Delta(R_i) < \Delta (R)$ for every $i$, we get $\Delta(L) \leq \Delta (R)$.
This and $R \subset L$ imply $\Delta (L) = \Delta (R)$.

Let $m_i, m_i'$ be the centers of the $(d-1)$-dimensional hemispheres $H_i$, $H_i'$ bounding $L_i$.
We maintain that at least one of these two centers, say $m_i$, belongs to the closure of $R \setminus R_i$.
The reason is that in the opposite case, there would be a neighborhood $N_i$ of $m_i$ such that $N_i \cap R_i = N_i \cap R$, which would imply that $H_i$ supports $R$ at $m_i$.
Moreover, $H_i'$ supports  $R$ at $m_i'$.
Hence $\Delta (R) = \Delta (L_i) = \Delta(R_i)$, in contradiction with  $\Delta (R_i) < \Delta(R)$. 

Since $m_i \in R \setminus R_i$ for $i=2, 3, \dots$,  we conclude that the sequence of points $m_2, m_3, \dots$ tends to $e$.
Consequently, $e$ is the center of a $(d-1)$-dimensional hemisphere bounding~$L$. 
\end{proof}

As noted in Remark 2 of \cite{L3}, besides the lune from Theorem \ref{e-reduced}, sometimes we have additional lunes $L' \supset R$ of thickness $\Delta (R)$ through $e$. 
Then $e$ is not in the middle of a $(d-1)$-dimensional hemisphere bounding $L'$. 
This happens, for instance, when $R$ is a spherical regular triangle $T_\alpha$ with $\alpha \leq {\pi \over 2}$.

By Theorem \ref{e-reduced} we obtain the following spherical analog of a theorem from the paper \cite{Gr} by Groemer (see also Corollary 1 in \cite{L0} and \cite{LaMa1}).
This theorem is presented in \cite{L5} and earlier for $d=2$ in \cite{LaMu1}.

\begin{thm}\label{smooth}  
Every reduced spherical body on $S^d$ of thickness greater than $\frac{\pi}{2}$ is smooth. 
\end{thm}

The thesis of Theorem \ref{smooth} does not hold for thickness at most $\frac{\pi}{2}$.
For instance, the regular spherical triangle of any thickness less than or equal to $\frac{\pi}{2}$ is reduced but not smooth.

\vskip0.21cm

\section{More about reduced bodies on $S^2$} \label{Reduced on $S^2$}

\vskip0.11cm
The notion of order of supporting hemispheres of a convex body $C \subset S^2$ presented in the next paragraph is needed for the later description of the boundaries of reduced bodies on $S^2$.
For this definition the classical notion of the polar set $F^\circ = \{ p : F\subset H(p)\}$ of a set $F\subset S^d$ is needed.
It is known that if $C$ is a spherical convex body, then $C^\circ$ is also a convex body.
Observe that $\bd (C^\circ) =\{ p: H(p)$ is a supporting hemisphere of $C\}$.

Let $C \subset S^2$ be a spherical convex body and $X=H(x)$, $Y=H(y)$, $Z=H(z)$ be different supporting hemispheres of $C$.
If $x,y,z$ are in this order on the boundary of $\bd (C^\circ)$ of $C$ when viewed from the inside of this polar, then (after \cite{LaMu1}) we write $\prec \!XYZ$ and say that $X, Y, Z$ {\bf support $C$ in this order}\index{order of hemispheres}.
Clearly, $\prec \! XYZ$, $\prec \! YZX$ and $\prec \! ZXY$ are equivalent.
The symbol $\preceq \! XYZ$ means that $\prec \!XYZ$ or $X=Y$ or $Y=Z$ or $Z=X$.
If $\prec XYZ$ (respectively, $\preceq XYZ$), then we say that {\bf $Y$ supports $C$ strictly between} (respectively {\bf between}) {\bf $X$  and $Z$}\index{supporting between}.

Let $X$ and $Z$ be hemispheres supporting $C$ at $p$ and let $\preceq \! XYZ$ for every hemisphere $Y$ supporting $C$ at $p$. Then $X$ is said to be the {\bf right supporting hemisphere at $p$}\index{right supporting hemisphere} and $Z$ is said to be the {\bf left supporting hemisphere at $p$}\index{left supporting hemisphere}. 
The left and right supporting hemispheres of $C$ are called {\bf extreme supporting hemispheres of $C$}\index{extreme supporting hemisphere}. 

Here are a lemma and a theorem from \cite{LaMu1}.

\begin{lem} \label{xxx}  
Let $C \subset S^2$ be a convex body with $\Delta (C) < \frac{\pi}{2}$ and let hemispheres $N_i$ with $\width_{N_i} = \Delta (C)$ for $i=1,2,3$ support $C$.
Then $\prec N_1N_2N_3$ if and only if $\prec N_1^*N_2^*N_2^*$.  
\end{lem}

\begin{thm} \label{main} 
Let $R \subset S^2$ be a reduced spherical body with $\Delta (R) < {\pi \over 2}$.
Let $M_1$ and $M_2$ be supporting hemispheres of $R$ such that ${\rm width}_{M_1} (R) = \Delta (R) = {\rm width}_{M_2} (R)$ and ${\rm width}_{M} (R) > \Delta (R)$ for every hemisphere $M$ satisfying $\prec \!M_1MM_2$.
Consider lunes $L_1 = M_1 \cap M_1^*$ and $L_2 = M_2 \cap M_2^*$.  
Then the arcs $a_1a_2$ and $b_1b_2$ are in $\bd (R)$, where $a_i$ is the center of $M_i/M_i^*$ and $b_i$ is the center of $M_i^*/M_i$ for $i=1,2$.
Moreover, $|a_1a_2| = |b_1b_2|$.
\end{thm}

The union of triangles $a_1a_2c$ and $b_1b_2c$, where $c$ denotes the intersection point of $a_1a_2$ and $b_1b_2$, is called a {\bf butterfly}\index{butterfly} and the pair $a_1a_2$, $b_1b_2$ is called {\bf a pair of arms}\index{pair of arms} of this butterfly (for instance see Example \ref{ExRedBodyOn$S^2$} with Figure \ref{fig:ExReduced}).
These names are analogous to the names from \cite{L2} for $E^2$, and for normed planes from \cite{FaLa} and \cite{LaMa2}.

Theorem \ref{main} allows us to describe the structure of the boundary of every reduced body on $S^2$ with $\Delta (R) < {\pi \over 2}$. 
Namely, we conclude that the boundary consists of pairs of arms of butterflies (in particular, the union of some arms of two butterflies may form a longer arc) and from some ``opposite" pieces of spherical curves of constant width (considered in Section \ref{Constant Width}).
We obtain these ``opposite" pieces always when ${\rm width}_M (R) = \Delta (R)$ for all $M$ fulfilling $\preceq \!M_1MM_2$, where $M_1$ and $M_2$ are two fixed supporting semi-circles of $R$.

The following proposition and two theorems are also taken from \cite{LaMu1}.

\begin{pro} \label{a1=a2} 
Let $R \subset S^2$ be a reduced body with $\Delta (R) < \frac {\pi}{2}$.
Assume that ${\rm width}_{M_1} (R) = \Delta (R) = {\rm width}_{M_2} (R)$, where $M_1$ and  $M_2$ are two fixed supporting hemispheres of $R$. 
Denote by $a_i$ and $b_i$ the centers of respectively the semi-circles $M_i/M_i^*$ and $M_i^*/M_i$ for $i=1,2$.
Assume that $a_1=a_2$ and that $\preceq \!M_1MM_2$ for every $M$ supporting $R$ at this point.
Then the shorter piece of the spherical circle with the center $a_1 = a_2$ and radius $\Delta (R)$ connecting $b_1$ and $b_2$ is in $\bd(R)$. 
Moreover, ${\rm width}_M (R) = \Delta (R)$ for all such $M$. 
\end{pro}

\begin{thm} \label{segment} 
Let $R$ be a reduced body with $\Delta (R) < \frac {\pi}{2}$.
Assume that $M$ is a supporting hemisphere of $R$ such that the intersection of $\bd(M)$ with $\bd (R)$ is a non-degenerate arc  $x_1x_2$. 
Then ${\rm width}_M (R) = \Delta (R)$, and the center of $M/M^*$ belongs to $x_1x_2$.
\end{thm}

\begin{thm} \label{extremehemisphere} 
 If $M$ is an extreme supporting hemisphere of a reduced spherical body $R \subset S^2$, then ${\rm width}_M (R) = \Delta (R)$. 
\end{thm}

After this theorem, the following problem is quite natural.

\begin{prob} 
Define the notion of extreme supporting hemispheres of a convex body $C \subset S^d$ and generalize Theorem \ref{extremehemisphere} to reduced convex bodies on $S^d$.
\end{prob}

For any maximum piece $\buildrel \frown \over {gh}$ of the boundary of a reduced body $R \subset S^2$ with $\Delta (R) < \frac{\pi}{2}$ which does not contain any arc the following claim from \cite{L9} holds true.

\begin{cla} \label{yyy}
If a hemisphere $K$ supports $R$ at a point of $\buildrel \frown \over {gh}$, 
then 
$\width_K (R) = \Delta(R)$.
\end{cla}

By Lemma \ref{xxx} and Claim \ref{yyy} we see that all the points at which our hemispheres $K^*$ touch $R$ form a  curve $\buildrel \frown \over {g'h'}$ in $\bd(R)$.
We call it the curve {\it opposite to the curve} $\buildrel \frown \over {gh}$.
Vice-versa, from this lemma we obtain that $\buildrel \frown \over {g'h'}$ determines $\buildrel \frown \over {gh}$.
So we say that $\buildrel \frown \over {gh}$ 
and $\buildrel \frown \over {g'h'}$ is a {\bf pair of opposite curves of constant width $\Delta(R)$}.  
As a recapitulation we obtain the following claim from \cite{L9}.
There the next theorem is also given.

\begin{cla}
The boundary of a reduced body $R \subset S^2$ of thickness below $\frac{\pi}{2}$ consists of countably many pairs of arms of butterflies and of countably many pairs of opposite pairs of curves of constant width $\Delta(R)$.
\end{cla}

\begin{thm} \label{approxReduced}
Let $R \subset S^2$ be a reduced body $R \subset S^2$ of thickness at most $\frac{\pi}{2}$.
For arbitrary $\varepsilon >0$ there exists a reduced body $R_\varepsilon \subset S^2$ with $\Delta (R_\varepsilon) = \Delta(R)$ whose boundary consists only of pairs of arms of butterflies and pieces of circles of radius $\Delta(R)$, such that the Hausdorff distance between $R_\varepsilon$ and $R$ is at most $\varepsilon$.
\end{thm}

Theorem \ref{e-reduced} leads to the following questions asked in \cite{L3} for $S^2$, which is rewritten here more generally for $S^d$.

\vskip0.21cm
\begin{prob} \label{lune passes}
Is it true that through every boundary point $p$ of a reduced body $R \subset S^2$ passes a lune $L \supset R$ of thickness $\Delta(R)$? 
\end{prob}

\vskip0.1cm
A consequence of a positive answer would be that every reduced body $R \subset S^2$ is an intersection of lunes of thickness $\Delta(R)$.

The example of the regular spherical triangle of thickness less than $\frac{\pi}{2}$ shows that the answer is negative when we additionally require that $p$ is the center of one of the two $(d-1)$-dimensional hemispheres bounding $L$. 

Theorem \ref{center} shows that for bodies of constant width on $S^d$ the answer to Problem \ref{lune passes} is positive.

After \cite{L1} recall that any reduced convex body $R \subset E^2$ is contained in a disk of radius $\frac{1}{2}\sqrt {\Delta (R)}$.
Its spherical analog is given in the main result of the paper \cite{Mus} by Musielak.

\begin{thm} \label{Th2Mus} 
Every reduced spherical body $R$ of thickness at most $\frac{\pi}{2}$ is contained in a disk of radius ${\rm arc tan} \big(\sqrt 2 \cdot \tan \frac{\Delta(R)}{2}\big)$.
\end{thm}

The proof of this theorem uses Theorem III from the paper \cite{Mo} by Moln\'ar concerning a spherical variant of Helly's theorem. This permits to avoid using the following consequence of our Corollary \ref{circumscribed}.

Recall that every reduced polygon $R \subset E^2$ is contained in a disk of radius $\frac{2}{3} \Delta(R)$ as shown in  Proposition from \cite{L1}.
This generates the subsequent problem for the two dimensional sphere.

\vskip0.2cm
\begin{prob}
What is the smallest radius of a disk which contains every reduced polygon of a given thickness on $S^2$?
\end{prob}

A general question is what properties of reduced bodies in $E^d$, and especially in $E^2$, (see \cite{L0} and \cite{LaMa1}) can be reformulated and proved for reduced spherical 
convex bodies.

\vskip0.21cm

\section{Spherical reduced polygons}\label{Polygons}

\vskip0.11cm
In this section we recall a number of facts on reduced polygons. 
But the following theorem and propositions from \cite{L4} are true for the more general situation of all spherical convex polygons.

\begin{thm} \label{Prop.2.1.[L3]} 
Let $V \subset S^2$ be a spherically convex polygon and let $L$ be a lune of thickness $\Delta (V)$ containing $V$.
Then at least one of the two semicircles bounding $L$ contains a side of $V$. 
If $\Delta (V) < {\pi \over 2}$, then the center of this semicircle belongs to the relative interior of this side.
\end{thm}

\begin{proof} 
Our lune $L$ has the form $G \cap H$, where $G$ and $H$ are different non-opposite hemispheres (see Figure \ref{fig:rotation}).

\vskip0.25cm
\begin{figure}[htbp]
\hskip1.95cm \includegraphics[width=9.265cm,height=4.335cm]{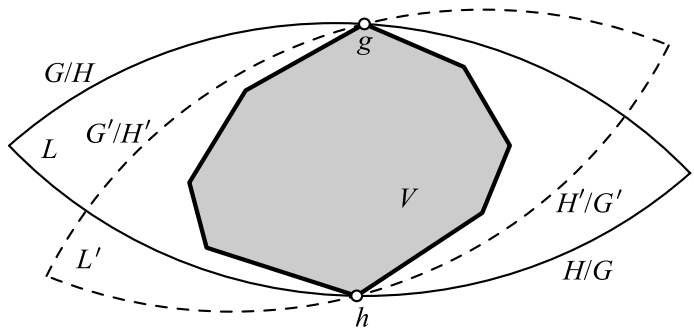} \\  
\caption{Rotation of $G$ around $g$ and rotation of $H$ around $h$} \label{fig:rotation}
\end{figure}

Let us prove the first statement.
Assume the contrary, i.e., that each of the two semicircles bounding $L$ contains exactly one vertex of $V$.
Denote by $g$ the vertex contained by $G/H$, and by $h$ the vertex contained by $H/G$. 
Apply Claim \ref{BothTheCenters}. 
We get that $g$ is the center of $G/H$ and $h$ is the center of $H/G$. 
Rotate $G$ around $g$, and $H$ around $h$, both according to the same orientation and by the same angle. 
Since $G/H$ and $H/G$ strictly support $V$, after any sufficiently small rotation of both with the same orientation, we obtain a new lune $L' = G' \cap H'$, where $G'$ and $H'$ are the images of $G$ and $H$ under corresponding rotations (see Figure 14). 
It still contains $V$ and has in common with the boundary of $V$ again only $g$ and $h$. 
But since now $g$ and $h$ are not the centers of the semicircles $G'/H'$, $H'/G'$ which bound $L'$ and since they are symmetric with respect to the center of $L'$, we get $\Delta (L') < \Delta (L)$. 
This contradicts the assumption $\Delta (L) = \Delta (V)$.
We see that the first part of our theorem holds true.

For the proof of the second part, assume the opposite, that is, that the center of the semicircle bounding $L$ does not belong to the relative interior of the side $S$ contained in this semicircle. 
This and Corollary \ref{Cor2LasWIDTH} imply that it must be an end-point of $S$, and thus a vertex of $V$. 
Apply $\Delta (V) < {\pi \over 2}$. 
Let us rotate $G$ around $g$, and $H$ around $h$, both with the same orientation, such that after sufficiently small rotations we get a lune narrower than $L$ still containing $V$.
This contradiction shows that our assumption at the beginning of this paragraph is false. 
\end{proof}

If $\Delta (V) = {\pi \over 2}$, then the center of the semicircle in the formulation of Theorem \ref{Prop.2.1.[L3]} may be an end-point of a side of $V$. 
This holds, for instance, for the regular triangle of thickness~${\pi \over 2}$.

Below we see a proposition needed for the proof of the subsequent theorem, both from \cite{L4}. 

\begin{pro} \label{Prop.2.2.[L3]} 
Let $V \subset S^2$ be a spherically convex polygon of thickness greater than $\pi \over 2$ and let $L$ be a lune of thickness $\Delta (V)$ containing $V$.
Then each of the two semicircles bounding $L$ contains a side of $V$.
Moreover, both centers of these semicircles belong to the relative interiors of the contained sides of~$V$.
\end{pro}

\begin{thm} \label{ReducedPolygonThickness} 
Every reduced spherical polygon is of thickness at most $\pi \over 2$.
\end{thm}

For a convex odd-gon $V = v_1v_2\dots v_n$ by the {\bf opposite side to the vertex $v_i$}\index{opposite side} we mean the side $v_{i + (n-1)/2}v_{i + (n+1)/2}$.
The indices are taken modulo~$n$.
 
The following fact is proved for any reduced polygon $V$ with $\Delta(V) < \frac{\pi}{2}$ in \cite{L4}. 
In Remark on p. 377 of \cite{L5}, applying Theorem \ref{equivalent}, this fact is transferred also for the case when $\Delta(V) = \frac{\pi}{2}$.
Let us add that a little later this fact for $\Delta(V) = \frac{\pi}{2}$ is given again in Theorem 3.1 of \cite{ChangLS}.
In summary, below we formulate this fact for every convex odd-gon on $S^2$.

\begin{thm} \label{ReducedPolygonProjection} 
A spherically convex odd-gon $V$ is reduced if and only if the projection of every vertex of $V$ on the great circle containing the opposite side belongs to the relative interior of this side and the distance of this vertex from this side is $\Delta (V)$. 
\end{thm}

\begin{figure}[htbp]
\hskip2.89cm \includegraphics[width=7.1cm,height=7.1cm]{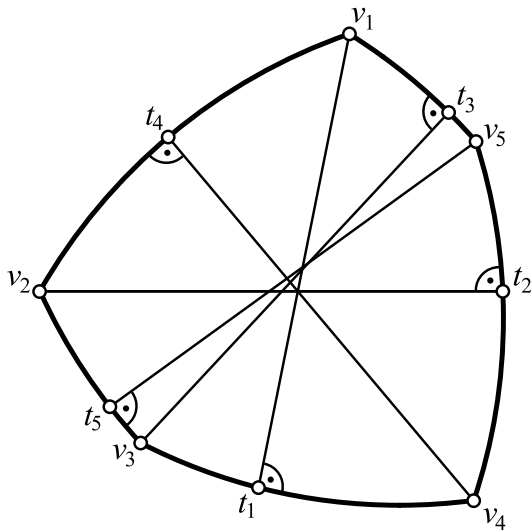} \\ 
\caption{A reduced spherical pentagon} \label{fig:pentagon}
\end{figure}

Figure \ref{fig:pentagon} shows a spherically convex pentagon. 
Theorem \ref{ReducedPolygonProjection} permits to recognize that it is reduced.

\begin{cor} \label{Cor3.3} 
Every spherical regular odd-gon of thickness at most ${\pi \over 2}$ is reduced.
\end{cor}

\begin{cor} \label{Cor3.4} 
The only reduced spherical triangles are the regular triangles of thickness at most~$\pi \over 2$. 
\end{cor}

\begin{cor} \label{Cor3.5} 
If $K$ is a supporting hemisphere of a reduced spherical polygon $V$ whose bounding circle contains a side of $V$, then ${\rm width}_K (V) = \Delta (V)$.
\end{cor}

The above three corollaries from our Theorem \ref{ReducedPolygonProjection} are given in \cite{L4}.
The following fact is also proved in \cite{L4} for $\Delta (V) < \frac{\pi}{2}$ and additionally for $\Delta (V) = \frac{\pi}{2}$ in the Remark of \cite{L5}. 
It is proved again in \cite{ChangLS}.

\begin{cor} \label{Cor3.6} 
For every reduced odd-gon $V= v_1v_2 \dots v_n$  with $\Delta (V) \leq {\pi \over 2}$ we have $|v_it_{i + (n+1)/2}| = |t_iv_{i + (n+1)/2}|$, for $i=1, \dots , n$, where $t_i$ denotes the projection of the vertex $v_i$ on the opposite side of $V$. 
\end{cor}

In addition, the paper \cite{ChangLS} by Chang, Liu and Su gives the following converse implication.

\begin{pro}
Every convex odd-gon $V= v_1v_2 \dots v_n$ on $S^2$ satisfying $|v_it_{i + (n+1)/2}| = |t_iv_{i + (n+1)/2}| = \frac{\pi}{2}$ for $i=1, \dots , n$ is reduced.
\end{pro}

\vskip0.15cm
Its proof is based on the fact given later in Proposition \ref{equivalent} that every reduced body of thickness $\frac{\pi}{2}$ is a body of constant width.
Thus this proof applying to bodies of constant width cannot be repeated for reduced bodies of thickness below $\frac{\pi}{2}$.
So let us formulate the following problem.

\vskip0.1cm
\begin{prob}
Consider an arbitrary convex odd-gon $V= v_1v_2 \dots v_n$ on $S^2$ of thickness below $\frac{\pi}{2}$ which satisfies $|v_it_{i + (n+1)/2}| = |t_iv_{i + (n+1)/2}|$ for every $i \in \{1, \dots , n\}$.  
Is $V$ a reduced polygon?
\end{prob}

Theorem \ref{ReducedPolygonThickness} and Corollary \ref{Cor3.6} together permit to construct reduced spherical polygons of arbitrary thickness at most $\pi \over 2$.
For instance see the reduced spherical pentagon in Figure \ref{fig:pentagon}. 

In the proof of Corollary \ref{Cor3.6} it is shown that the triangles $v_io_it_{i + (n+1)/2}$ and $v_{i + (n+1)/2}o_it_i$, where $o_i$ is the common point of the arcs $v_it_i$ and $v_{i+(n+1)/2}t_{i+(n+1)/2}$, are symmetric with respect to a great circle. 
This implies the following corollary.

\begin{cor} \label{Cor3.7} 
In every reduced odd-gon $V= v_1v_2 \dots v_n$, for every $i \in \{1, \dots , n \}$ we have $\angle t_{i + (n+1)/2}v_it_i = \angle t_iv_{i + (n+1)/2}t_{i + (n+1)/2}$. 
\end{cor}

Also the following three corollaries are obtained in \cite{L4} for $\Delta (V) < {\pi \over 2}$ and in Remark \cite{L5} expanded also for $\Delta (V) = {\pi \over 2}$. 

By Corollary \ref{Cor3.6} applied $n$ times, the sum of the lengths of the boundary spherical arcs of $V$ from $v_i$ to $t_i$ (with the positive orientation) is equal to the sum of the boundary spherical arcs of $V$ from $t_i$ to $v_i$ (with the positive orientation).  
Let us formulate this statement as the following corollary. 

\begin{cor} \label{Cor3.8} 
Let $V$ be a reduced spherical $n$-gon and let $i \in \{1,\dots , n\}$. 
For every $i \in \{1,\dots , n\}$ the spherical arc $v_it_i$ divides in half the perimeter of~$V$.
\end{cor}

\begin{cor} \label{Cor3.10}   
For every reduced polygon and every $i$
we have $\alpha_i > {\pi \over 6} +E$ and  $\beta_i < {\pi \over 6} +E$, where $E$ denotes the excess of the triangle $v_it_iv_{i + (n+1)/2}$.
\end{cor} 

In order to formulate the next corollary, for every $i \in \{1, \dots , n\}$ let us put $\alpha_i = \angle v_{i+1}v_it_i$ and $\beta_i = \angle t_iv_iv_{i + (n+1)/2}$ in a spherically convex odd-gon $V= v_1v_2 \dots v_n$. 

The angle $\alpha_i$ is greater, but the one for $\beta_i$ is smaller than these given in Theorem 8 of \cite{L0} for a reduced $n$-gon in $E^2$.
For instance, in the regular spherical triangle of sides and thickness $\pi/2$ we have $\alpha_i = \beta_i = \pi/4$.

\begin{cor} \label{Cor3.9}  
If $V= v_1v_2 \dots v_n$ is a reduced spherical polygon with $\Delta (V) < {\pi \over 2}$, then $\beta_i \leq \alpha_i$ for every $i \in \{1, \dots , n\}$.
\end{cor}

Example 3.3 from \cite{ChangLS} constructs a spherically convex pentagon of thickness $\frac{\pi}{2}$ answering the question from \cite{L4} about the existence of non-regular reduced polygons of thickness $\frac{\pi}{2}$. 

The paper \cite{L4} presents two conjectures on the perimeter of a reduced spherical polygon.
The first says that the perimeter of every reduced spherical polygon $V$ is not larger than that of the spherical regular triangle $T$ of the same thickness, so at most $6 \; {\rm arc cos} {\cos \Delta (T) + \sqrt {8+ \cos^2 \Delta (T)} \over 4}$, and that it is attained only for this regular triangle. 
The second says that from amongst all reduced spherical polygons of fixed thickness and with at most $n$ vertices only the regular spherical $n$-gon has  minimal perimeter.
For the special case when $\Delta(V) = \frac{\pi}{2}$ both conjectures are confirmed by Chang, Liu and Su \cite{ChangLS} in the following two theorems. 
 
\begin{thm} \label{CLS 4.1}  
The perimeter of every reduced spherical polygon of thickness $\frac{\pi}{2}$
is less than or equal to that of the regular spherical triangle of thickness $\frac{\pi}{2}$
and the maximum value is attained only for this regular spherical triangle. 
\end{thm}

\begin{thm} \label{CLS 4.3}  
The regular spherical $n$-gon of thickness $\frac{\pi}{2}$ has minimum perimeter among all reduced spherical $k$-gons of thickness $\frac{\pi}{2}$. 
\end{thm}

Since the second conjecture is not proved yet for thickness below $\frac{\pi}{2}$, we propose the next problem.

\begin{prob}
Is it true that the perimeter of every reduced spherical polygon of thickness below $\frac{\pi}{2}$ is less than or equal to that of the regular spherical triangle of the same thickness? 
Is the maximum value attained only for this regular spherical triangle?
Is it true that the regular spherical $n$-gon of any fixed thickness below $\frac{\pi}{2}$ has minimum perimeter among all reduced spherical $n$-gons of the same thickness?
\end{prob}

The following two facts are also presented by Chang, Liu and Su \cite{ChangLS}.

\begin{cla} \label{CLS 4.2}  
The regular spherical $n$-gon of thickness $\frac{\pi}{2}$ has minimum perimeter among all regular spherical $n$-gons of thickness $\frac{\pi}{2}$.
\end{cla}

\begin{cor} \label{CLS 4.4} 
The regular spherical $n$-gon of thickness $\frac{\pi}{2}$ has minimum perimeter among all reduced spherical $k$-gons of thickness $\frac{\pi}{2}$, where $k$ and $n$ are odd natural numbers such that $3 \leq k \leq n$. 
\end{cor}

Finally in this section, let us recall a few facts about the area of spherical convex polygons and especially if they are reduced.

A generalization of the classical formula of Girard for the area of a spherical triangle is presented by Todhunter in Part 99 of \cite{Tod}. 
In particular, he  says that the area of an arbitrary convex $n$-gon on $S^2$ is equal to the sum of its internal angles diminished by $(n - 2)\pi$.

The following theorem is proved for polygons of thickness $\frac{\pi}{2}$ by Chang, Liu and Su in \cite{ChangLS}, while for the polygons of any thickness smaller than $\frac{\pi}{2}$ by Liu, Chang and Su in \cite{LiuChangSu}.
 
\begin{thm} \label{CLS 5.2} 
The regular spherical $n$-gon has maximum area amongst all regular spherical $k$-gons with odd numbers $k, n$ and $3 \leq k \leq n$.
\end{thm}

The paper \cite{LiuChangSu} by Liu, Chang and Su confirms two conjectures from \cite{L4} by proving the following two facts.

\begin{thm} \label{CLS 5.3}
Every reduced spherical non-regular $n$-gon of any fixed thickness has area smaller than the regular spherical $n$-gon of this thickness.
\end{thm}

\begin{cor} \label{CLS 5.4} 
The area of arbitrary reduced spherical polygon $V$ is smaller than $2(1 - \cos \frac{\Delta(V)}{2})\pi$.
\end{cor}

\vskip0.21cm

\section{Diameter of convex bodies and reduced bodies} \label{Diameter}

\vskip0.11cm
The following theorem and proposition from \cite{L3} are spherical analogs of the classical facts in the Euclidean space.

\begin{thm} \label{MaxWidth=Diam} 
Let ${\rm  diam} (C) \leq {\pi \over 2}$ for a convex body $C \subset S^d$. 
We have 

\vskip-0.42cm
$$\max \{ {\rm width}_K (C); K {\rm \ is \ a \ supporting \ hemisphere \ of} \ C \} = {\rm  diam} (C).$$
\end{thm}

\vskip0.15cm
The following example from \cite{L3} shows that Theorem \ref{MaxWidth=Diam} does not hold without the
assumption ${\rm  diam} (C) \leq {\pi \over 2}$.

\vskip0.25cm
\noindent
{\bf Example 3.}
Let $T$ be an isosceles triangle with base of length $\lambda$ close to $0$ and the altitude perpendicular to it of length $\mu \in ({\pi \over 2}, \pi)$ (see Figure 7). 
Denote by $w$ the center of the base, by $v$ the opposite vertex of $T$ and by $y, z$ the remaining vertices of $T$. 
Claim \ref{DistInLunes} implies that $wv$ is the diametrical arc of $T$.
Take the hemisphere $K$ orthogonal to $vw$ supporting $T$ at $w$.  
Denote by $k$ the center of $K$.
Clearly, $k \in wv$, so $k$ is in the interior of $T$.
Let $\rho$ be the radius of the largest disk $B$ with center $k$ contained in $T$.
The radius $\rho$ of $B$ is arbitrarily close to $0$, as $\lambda$ is sufficiently small.
Applying Part III of Theorem \ref{I-III} we conclude that the width of $T$ determined by $K$ is ${\pi \over 2} + \rho$  which is realized for the hemisphere $K^*$ whose boundary contains $vy$ (that is, $\width_K = \angle vyz$). 
Hence it may be arbitrarily close to $\pi \over 2$, as $\lambda$ is sufficiently small.
On the other hand, the diameter $|wv|$ of $T$ may be arbitrarily close to $\pi$, as $\mu$ is sufficiently close to $\pi$.

\vskip-0.23cm
\begin{figure}[htbp]
\hskip0.8cm \includegraphics[width=11.28cm,height=4.24cm]{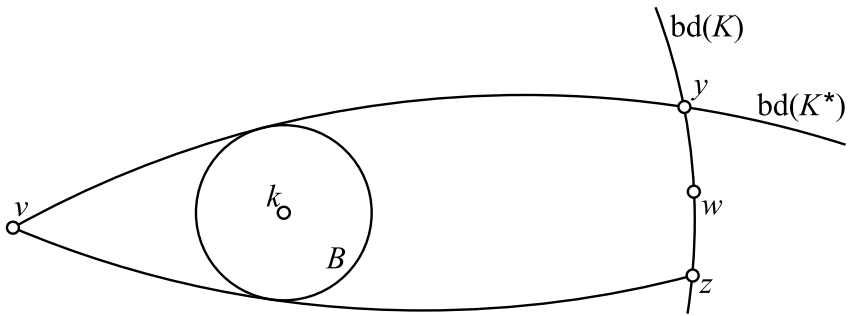} \\ 
\caption{Illustration to Example 3}  
\end{figure}

\vskip0.4cm
\begin{pro} \label{Prop1L2}
Let ${\rm  diam} (C) > {\pi \over 2}$ for a convex body $C \subset S^d$.  
We have 

\vskip-0.4cm
$$\max \{ {\rm width}_K (C); K {\rm \ is \ a \ supporting \ hemisphere \ of} \ C \} \leq {\rm  diam} (C).$$
\end{pro}

The paper \cite{L5} presents the following corollary from this proposition.

\begin{cor} \label{DeltaATmostDiam} 
For an arbitrary convex body $C \subset S^d$ we have $\Delta (C) \le \diam (C)$. 
\end{cor}

Let us apply Claim \ref{diameter} for a convex body $C$ of diameter larger than $\frac{\pi}{2}$. 
Having in mind that the center $k$ of $K$ is in $pq$ and thus in $C$, by Part III of Theorem \ref{I-III} we obtain $\Delta (K \cap K^*) > \frac{\pi}{2}$. 
This gives the following corollary from \cite{LaMu2}, which implies the second one.

\begin{cor} \label{ortho} 
Let $C \subset S^d$ be a convex body of diameter larger than $\frac{\pi}{2}$ and let 
$\diam (C)$ be realized for points $p, q \in C$.
Take the hemisphere $K$ orthogonal to $pq$ at $p$ which supports $C$.
Then ${\rm width}_K (C) > \frac{\pi}{2}$. 
\end{cor}

\begin{cor} \label{family} 
Let $C \subset S^d$ be a convex body of diameter larger than $\frac{\pi}{2}$ and let $\mathcal K$ denote the family of all hemispheres supporting $C$.  
Then we have $\max_{K \in \mathcal K} {\rm width}_K (C) > \frac{\pi}{2}$. 
\end{cor}

The set of extreme points of $C$ is denoted by $E(C)$.
Here is a property established in \cite{LaMu3}. 

\begin{cla} \label{diamE} 
For every convex body $C \subset S^2$ of diameter at most $\frac{\pi}{2}$ we have $\diam (E(C)) = \diam(C)$.
\end{cla}

The assumption that $\diam (C) \leq \frac{\pi}{2}$ is essential in Claim \ref{diamE}, as it follows from the example of a regular triangle of any diameter greater than $\frac{\pi}{2}$. 
The weaker assumption that $\Delta(C) \leq \frac{\pi}{2}$ is not sufficient, 
which we see taking in the part of $C$ any isosceles triangle $T$ with $\Delta(T) \leq \frac{\pi}{2}$ and legs longer than $\frac{\pi}{2}$ (so with base shorter than $\frac{\pi}{2}$). 
The diameter of $T$ is equal to the distance between the midpoint of the base and the opposite vertex of $T$. 
So $\diam (T)$ is greater than the length of each of the sides.

The above claim was established as a lemma needed in Theorem \ref{ineq} for $S^2$ only, and thus it concerns only bodies on $S^2$.
We expect that the $d$-dimensional variant of Claim  \ref{diamE} holds as well.

From Theorem \ref{smooth} and Corollary \ref{family} we get obtain following three corollaries for reduced bodies on $S^d$ (not only for bodies of constant width, as formulated in~\cite{LaMu2}).

\begin{cor} \label{diamover}
If $\diam (R) >  \frac{\pi}{2}$ for a reduced body of $R \subset S^d$, then $\Delta(R) >  \frac{\pi}{2}$.
\end{cor}

\begin{cor} \label{diamless}
For every reduced body of $R \subset S^d$ of thickness at most $\frac{\pi}{2}$ we have $\diam(R) \leq \frac{\pi}{2}$.
\end{cor}

\begin{cor} \label{H(p)}
Let $p$ be a point of a reduced body $R \subset S^d$ of thickness at most $\frac{\pi}{2}$. 
Then $R \subset H(p)$.
\end{cor}

By the way, the last corollary was shown earlier for $S^2$ in \cite{LaMu1}.

Note the obvious fact that the diameter of a convex body $C \subset S^d$ is realized only for some pairs of points of $\bd(C)$.
In the following four claims we recall some facts from \cite{L3}, \cite{L5} and \cite{LaMu2}.
Two first two concern the realization of the diameter of any spherical convex body.

\begin{cla}\label{diameter} 
Assume that the diameter of a convex body $C \subset S^d$ is realized for points $p$ and $q$. 
The hemisphere $K$ orthogonal to $pq$ at $p$ and containing $q \in K$ supports $C$.
\end{cla}

\begin{cla} \label{extreme20}  
Let $C \subset S^d$ be a convex body. 
If $\diam (C) < \frac{\pi}{2}$, then every two points of $C$ at the distance $\diam (C)$ are extreme.
If $\diam (C) = \frac{\pi}{2}$, then there are two points of $C$ at distance $\frac{\pi}{2}$ such that at least one of them is an extreme point of $C$. 
\end{cla}

Clearly, the first thesis in not true for $\diam (C) = \frac{\pi}{2}$.
After the proof of this claim in \cite{L5} the following amplification is also shown: {\it every convex body $C \subset S^2$ of diameter $\frac{\pi}{2}$ contains a pair of extreme points distant by $\frac{\pi}{2}$}.
Now, we also observe that this stronger statement is also true for $C \subset S^d$, as formulated in next claim.

\begin{cla} \label{extreme20+}
Every convex body $C \subset S^d$ of diameter $\frac{\pi}{2}$ contains a pair of extreme points distant by $\frac{\pi}{2}$.
\end{cla}

\begin{proof}
Consider a point $a \in C$ and an extreme point $b \in C$ guaranteed by the second part of Claim \ref{extreme20}.
If $a$ is also extreme, the thesis is true.
In the contrary case take the $(d-1)$-dimensional hemisphere $K$ bounding the lune $L$ from Claim \ref{C-subset-L} whose center is $a$.
By Claim \ref{extreme} there are at most $d$ extreme points of $C \cap K$ whose convex hull contains $a$.
Let $a'$ be one of them.
By the second part of Lemma \ref{DistInLunes}, we have $|a'b| = \frac{\pi}{2}$ as required.
\end{proof}

The following theorem from \cite{LaMu3} is analogous to the first part of Theorem 9 from \cite{L0} and confirms the conjecture from the bottom of p. 214 of \cite{L4}.

\begin{thm} \label{ineq} 
For every reduced spherical body $R \subset S^2$ with $\Delta (R) < \frac{\pi}{2}$ we have 

\vskip-0.4cm
$$\diam(R) \leq {\rm arccos} (\cos^2 \Delta (R)).$$
This value is attained if and only if $R$ is the quarter of a disk of radius $\Delta(R)$.
If $\Delta (R) \geq \frac{\pi}{2}$, then $\diam(R) = \Delta (R)$.
\end{thm}

From L'Hospital's rule we conclude that if $R$ is just a quarter of disk, then the ratio $[{\rm arc cos}(\cos^2 \Delta (R))]/\Delta (R)$ tends to $\sqrt 2$ as $\Delta (R)$ tends to $0$. 
Consequently, the limit factor $\sqrt 2$ is like in the planar case in Theorem 9 of \cite{L0}.

By the way, the weaker estimate $\diam(R) \leq 2 \arctan \left(\sqrt 2 \tan \frac{\Delta(R)}{2}\right)$ than the one in Theorem \ref{ineq} is a consequence of Theorem \ref{Th2Mus}.

Theorem \ref{ineq} implies the following proposition proved also in \cite{LaMu3}.

\begin{pro} \label{precise} 
Let $R \subset S^2$ be a reduced body. 
Then $\diam (R) < {\pi \over 2}$ if and only if $\Delta (R) < {\pi \over 2}$. 
Moreover, $\diam (R) = {\pi \over 2}$ if and only if $\Delta (R) = {\pi \over 2}$. 
\end{pro}

This proposition gives a more precise statement than the one $\diam (R) \leq \frac{\pi}{2}$ following from 
Corollary \ref {H(p)}.

An example from \cite{L1} shows that in $E^d$, where $d\geq 3$, there are reduced bodies of thickness $1$ and arbitrary large diameter. 
Below we formulate an analogous question for reduced bodies on $S^d$.

\vskip0.15cm
\begin{prob}
Is it true that for any number $q > 1$ there exists a reduced body on $S^3$ with the quotient of its diameter to its thickness at least $q$? 
Construct, if possible, such an example on $S^d$, with any $d\geq 3$.
\end{prob}

The next theorem is proved in \cite{L5}.

\begin{thm} \label{Delta-DIAM} 
For every reduced body $R\subset S^d$ such that $\Delta(R) \leq \frac{\pi}{2}$ we have $\diam (R) \leq \frac{\pi}{2}$.
Moreover,  if $\Delta(R) < \frac{\pi}{2}$, then $\diam (R) < \frac{\pi}{2}$.
\end{thm}

The special case of the first assertion of this theorem for $d=2$ is stated in the observation just before Proposition~1 of \cite{LaMu3}.

By the way, for $d=2$ the first assertion of our Proposition \ref{precise} is a special case of the first statement of Theorem \ref{iff} and the second assertion appears to be a consequence of just Theorem \ref{Delta-DIAM}.
Let us add that the approach of the proof of Theorem \ref{Delta-DIAM} is different from the argumentation in the proof of Proposition \ref{precise} which applies Theorem \ref{ineq} proved only for $d=2$.

From Theorems \ref{diam=w} and \ref{Delta-CW} the following proposition and corollary from \cite{L5} result.

\begin{pro} \label{diam=Delta} 
For every reduced body $R \subset S^d$ satisfying $\Delta (R) \geq \frac{\pi}{2}$ we have $\Delta (R) = \diam (R)$. 
\end{pro}

Observe that this proposition does not hold without the assumption that the body is reduced.

\begin{cor} 
If $\Delta (R) < \diam (R)$ for a reduced body $R  \subset S^d$, then both numbers are below $\frac{\pi}{2}$.
Moreover, $R$ is not a body of constant width. 
\end{cor}

Here are five equivalences from \cite{L5}.

\begin{thm}\label{iff} 
Let $R \subset S^d$ be a reduced body. 
Then
\vskip 0.1cm

{\rm (a)} $\Delta(R) =\frac{\pi}{2}$  if and only if $\diam (R) =\frac{\pi}{2}$,

{\rm (b)} $\Delta(R) \geq \frac{\pi}{2}$  if and only if $\diam (R) \geq \frac{\pi}{2}$,

{\rm (c)} $\Delta(R) > \frac{\pi}{2}$  if and only if $\diam (R) > \frac{\pi}{2}$,

{\rm (d)} $\Delta(R) \le \frac{\pi}{2}$  if and only if $\diam (R) \le \frac{\pi}{2}$, 

{\rm (e)} $\Delta(R) < \frac{\pi}{2}$  if and only if $\diam (R) < \frac{\pi}{2}$. 
\end{thm}

Let us add that (a) and (e) were established earlier for $d=2$ in \cite{LaMu3}.

As in Section  4 of \cite{LaMu2}, we say that a convex body $D \subset S^d$ of diameter $\delta$ is {\bf of constant diameter}\index{constant diameter} $\delta$ provided for arbitrary $p \in \bd (D)$ there exists $p' \in \bd (D)$ such that $|pp'| = \delta$ (more generally, this notion makes sense for a closed set $D$ of diameter $\delta$ in a metric space $M$, such that for all $x, z \in D$ and $y \in M$ with $|xy| + |yz| = |xz|$ we have $y \in D$, so for instance when $M$ is a Riemannian manifold). 
This is an analog of the notion of a body of constant diameter in Euclidean space considered by Reidemeister \cite{Re}.    

We get some bodies of constant diameter on $S^2$ as particular cases of the later Example \ref{prolong} by taking there any non-negative $\kappa < \frac{\pi}{2}$ and $\sigma = \frac{\pi}{4} - \frac{\kappa}{2}$. 
The following example from \cite{L5} presents a wider class of bodies of constant diameter $\frac{\pi}{2}$.

\begin{exa} 
{\rm Take a triangle $v_1v_2v_3 \subset S^2$ of diameter at most $\frac{\pi}{2}$. 
Put $\kappa_{12} =~|v_1v_2|, \kappa_{23}=|v_2v_3|, \kappa_{31}=|v_3v_1|$,
$\sigma_1 = \frac{\pi}{4} - \frac{\kappa_{12}}{2} + \frac{\kappa_{23}}{2} - \frac{\kappa_{31}}{2}$,
$\sigma_2 = \frac{\pi}{4} - \frac{\kappa_{12}}{2} - \frac{\kappa_{23}}{2} + \frac{\kappa_{31}}{2}$,
$\sigma_3 = \frac{\pi}{4} + \frac{\kappa_{12}}{2} - \frac{\kappa_{23}}{2} - \frac{\kappa_{31}}{2}$.
Here we consent only to triangles with the sum of  lengths of two shortest sides at most the length of the longest side plus $\frac{\pi}{2}$ (equivalently: with $\sigma_1 \geq 0$, $\sigma_2 \geq 0$ and $\sigma_3 \geq 0$).
Extend the following $v_1v_2$ up to $w_{12}w_{21}$ with $v_1 \in w_{12}v_2$, 
$v_2v_3$ up to $w_{23}w_{32}$ with $v_2 \in w_{23}v_3$, 
$v_3v_1$ up to sides: $w_{31}w_{13}$ with $v_3 \in w_{31}v_1$  
such that 
$|v_1w_{12}| = |v_1w_{13}| = \sigma_1$,
$|v_2w_{21}| = |v_2w_{23}| = \sigma_2$,
$|v_3w_{31}| = |v_3w_{32}| = \sigma_3$ (see Figure \ref{fig:ExConstantDiameter}).
Draw six pieces of circles: 
with center $v_1$ of radius $\sigma_1$ from $w_{12}$ to $w_{13}$ and of radius $\frac{\pi}{2} - \sigma_1$ from $w_{21}$ to $w_{31}$, 
with center $v_2$ of radius $\sigma_2$ from $w_{23}$ to $w_{21}$ and of radius $\frac{\pi}{2} - \sigma_2$ from $w_{32}$ to $w_{12}$, 
with center $v_3$ of radius $\sigma_3$ from $w_{31}$ to $v_{32}$ and of radius $\frac{\pi}{2} - \sigma_3$ from $w_{13}$ to $w_{23}$.
Clearly, the convex hull of these six pieces of circles is a body of constant diameter~$\frac{\pi}{2}$.

\vskip0.2cm
\begin{figure}[htbp]
\hskip2.75cm \includegraphics[width=2.9in]{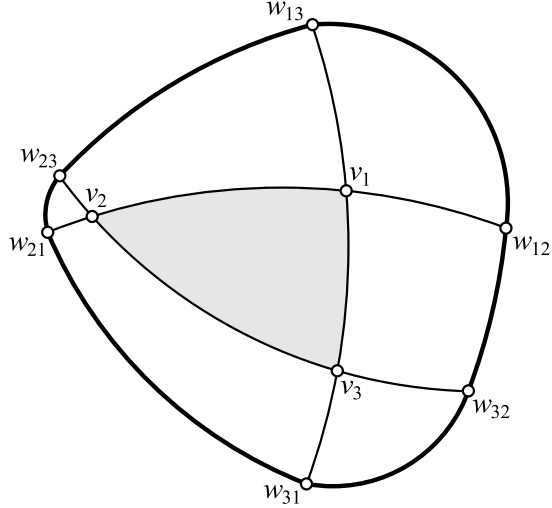} \\ 
\vskip-0.25cm
\caption{A spherical body of constant diameter}  \label{fig:ExConstantDiameter}
\end{figure}

Generalizing, take a convex odd-gon $v_1 \dots v_n \subset S^2$ of diameter at most $\frac{\pi}{2}$.
 For $i=1, \dots ,n$ put $\kappa_{i \ {i + (n-1)/2}}\ = |v_iv_{i + (n-1)/2}|$ (here and later we mean indices modulo $n$). 
Let $\sigma_i = \frac{\pi}{4} + \Sigma_{i=1}^n s_i\kappa_{i \ {i + (n-1)/2}}$, where $s_i = \frac{1}{2}$ if $|\frac{n+1}{2} -i| \leq \frac {n-3}{4}$ for $n$ of the form $3+4k$ and $|\frac{n}{2} - i| \leq \frac {n-1}{4}$ for $n$ of the form $5+4k$ (where $k=0,1,2, \dots$), and $s_i = -\frac{1}{2}$ in the opposite case.
We agree only for odd-gons with $\sigma_i \geq 0$ for $i=1, \dots , n$.
Prolong each diagonal $v_iv_{i + (n-1)/2}$ up to the arc $w_{i \ i + (n-1)/2} w_{i + (n-1)/2  \ i}$ such that $v_i \in w_{i \ i + (n-1)/2}v_{i + (n-1)/2}$ and $|v_iw_{i \ i + (n-1)/2}| = |v_iw_{i \ i + (n+1)/2}|= \sigma_i$.
For $i= 1, \dots , n$ we draw the piece of the circle with center $v_i$ of radius $\sigma_i$ from $w_{i \ i + (n-1)/2}$ to $w_{i \ i + (n+1)/2}$ and the piece of circle of radius $\frac {\pi}{2} - \sigma_i$ from $w_{i + (n-1)/2  \ i}$ to $w_{i + (n+1)/2  \ i}$. 
The convex hull of the union of these $2n$ pieces is a body of constant diameter $\frac{\pi}{2}$.}
\end{exa}

The proposition below from \cite{L6} presents a property of bodies of constant diameter, and consequently also for bodies of constant width and complete considered in the two next Sections.
It does not hold for $S^d$ with $d\geq 3$.
A not very complicated example of such a body of constant diameter on $S^3$ is the body obtained by the rotation of a quarter of a disk on $S^2$ around its axis of symmetry.  

\begin{pro} \label{intersect} 
Let $D \subset S^2$ be  a body of constant diameter. 
 Then every two diametral chords of $D$ intersect. 
\end{pro}

Note that Proposition \ref{realized} and Theorem \ref {DiamPolygon} below give some information on the diameter of reduced polygons.

The paper \cite{L0} addresses the problem of whether there exist reduced polytopes in $E^d$ for $d\geq 3$. 
It is positively answered by Gonzalez Merino, Jahn, Polyanskii and Wachsmuth in \cite{GJPW}.
For $S^d$ the analogous question has a specific positive answer.
The spherical simplex being the $1\over 2^d$-part of the ball of radius $\pi \over 2$ is a reduced spherical polytope (by the way, in $E^d$ with $d\geq 3$ all simplices are not reduced, as shown by Martini and Swanepoel \cite{MS}).  
This follows from the fact that it is a spherical body of constant width.
If we disregard this example, we still have the following problem. 

\begin{prob}
Do there exist spherical $d$-dimensional polytopes (possibly some simplices?) on $S^d$, where $d \geq 3$, different from the $1\over 2^d$-part of $S^d$? 
\end{prob}

Let us recall a proposition and a theorem from \cite{L4}.

\begin{pro} \label{realized} 
The diameter of any reduced spherical $n$-gon is realized only for some pairs of vertices whose indices (modulo $n$) differ by ${n-1}\over 2$ or ${n+1}\over 2$. 
\end{pro}

Here is a sharper estimate for reduced polygons in place of reduced bodies on $S^2$ than the estimate from Theorem \ref{ineq}. 

\begin{thm} \label{DiamPolygon} 
For every reduced spherical polygon $V \subset S^2$ we have

$$\diam(V) \leq {\rm arc cos} \Big(\sqrt {1 - {\sqrt 2 \over 2} \sin \Delta (V)} \cdot \cos \Delta (V) \Big),$$

\noindent 
{\it with equality for the regular spherical triangle in the part of $V$.}
\end{thm}

For instance, from this theorem we see that every reduced spherical polygon of thickness $\pi \over 4$ is of diameter at most $\pi \over 3$ and that the regular triangle of thickness $\pi \over 4$ has diameter~$\pi \over 3$.

By the way, from Theorem \ref{DiamPolygon} we get a slightly better estimate of the radius of a disk
covering a reduced spherical polygon than the one in Theorem \ref{Th2Mus}. 
The resulting formula is very complicated, so we omit it here.
Let us give only an example; for $\Delta (R) = 60^{\circ}$ it gives the estimate about $(38,38)^\circ$ for the radius, while the estimate from Theorem \ref{Th2Mus} gives about  $(39,23)^\circ$.

\begin{prob}
We conjecture that the equality in Theorem \ref{DiamPolygon} holds only for regular triangles.
\end{prob}

By the way, recently Horv\'ath \cite{Hor} showed that every body of constant diameter in the hyperbolic plane is a body of constant width.

\vskip0.21cm

\section{Bodies of constant width} \label{Constant Width}

\vskip0.11cm
We say that a convex body $W \subset S^d$ is of {\bf constant width}\index{constant width} $w$ provided that for every supporting hemisphere $K$ of $W$ we have $\width_K (W) = w$.

This notion is analogous to the classical notion of bodies of constant width in $E^n$. 
The oldest paper on this subject we found is the famous paper of Barbier \cite{Barb} which proves that the area of all planar bodies of the same constant width is the same.
The book \cite{MMO} by Martini, Montejano and Oliveros gives a wide survey of results on bodies of constant width in various structures.
Also Section 11 of the book \cite{Lay} by Lay and
Part 2.5 of the book \cite{Toth} by Toth, besides classical facts on bodies of constant width recall a number of less known properties and exercises of them.
Many properties of bodies of constant width in the plane are recalled in the popular book \cite{YaBo} by Yaglom and Boltyanskij.

In particular, spherical balls of radius smaller than $\frac{\pi}{2}$ are spherical bodies of constant width. 
Also every spherical Reuleaux odd-gon is a convex body of constant width.
Each of the ${1 \over 2^d}$-th part of any ball on $S^d$ is a spherical body of constant width $\frac{\pi}{2}$, which easily follows from the definition of a body of constant width.   

After \cite{LaMu2} we remind the reader an example of a spherical body of constant width on the sphere $S^3$.

\begin{exa} \label{RedBodyS3} 
{Take a circle $X \subset S^3$ (i.e., a set congruent to a circle in $E^2$) of a positive diameter $\kappa < \frac{\pi}{2}$, and a point $y \in S^3$ at the distance $\kappa$ from every point $x \in X$. 
Prolong every spherical arc $yx$ by a distance $\sigma \leq {\pi \over 2} - \kappa$ up to points $a$ and $b$ so that $a, y, x, b$ are on one great circle in this order.
All these points $a$ form a circle $A$, and all points $b$ form a circle $B$. 
On the sphere on $S^3$ of radius $\sigma$ whose center is $y$ take the ``smaller" part $A^+$ bounded by the circle $A$. 
On the sphere on $S^3$ of radius $\kappa + \sigma$ with center $y$ take the ``smaller" part $B^+$ bounded by $B$.
For every $x \in X$ denote by $x'$ the point on $X$ such that $|xx'| = \kappa$.
Prolong every $xx'$ up to points $d, d'$ so that $d, x, x', d'$ are in this order and $|dx|= \sigma =|x'd'|$.
For every $x$ provide the shorter piece $C_x$ of the circle with center $x$ and radius $\sigma$ connecting the $b$ and $d$ determined by $x$ and also the shorter piece $D_x$ of the circle with center $x$ of radius $\kappa +\sigma$ connecting the $a$ and $d'$ determined by $x$. 
Denote by $W$ the convex hull of the union of $A^+$, $B^+$ and all pieces $C_x$ and $D_x$. 
This is a body of constant width $\kappa + 2\sigma$ (hint: for every hemisphere $H$ supporting $W$ and every $H^*$ the centers of $H/H^*$ and $H^*/H$ belong to $\bd (W)$ and the arc connecting them passes through one of our points $x$, or through the point $y$).}
\end{exa}

Recall the following related question from p. 563 of \cite{L3}.

\begin{prob}
Is a convex body $W \subset S^d$ of constant width provided every supporting hemisphere $G$ of $W$ determines at least one hemisphere $H$ supporting $W$ such that $G \cap H$ is a lune with the centers of $G/H$ and $H/G$ in $\bd (W)$? 
\end{prob}

Let us formulate a comment from \cite{L3}, p. 562--563 as the following claim.

\begin{cla} \label{determines hemisphere}
 If $W\subset S^d$ is a body of constant width, then every supporting hemisphere $G$ of $W$ determines a supporting hemisphere $H$ of $C$ for which $G\cap H$ is a lune such that the centers of $G/H$ and $H/G$ belong to the boundary of $C$. 
\end{cla}

\begin{prob} \label{opposite}
Is the contrary true?
More precisely, is a convex body $C \subset S^d$ of constant width provided every supporting hemisphere $G$ of $C$ determines at least one hemisphere $H$ supporting $C$ such that $G \cap H$ is a lune with the centers of $G/H$ and $H/G$ in the boundary of $C$?
\end{prob}

The above claim implies the next one from \cite{L3}.

\begin{cla} \label{lunesCW}
Every spherical body $W$ of constant width is an intersection of lunes of thickness $\Delta (W)$ such that the centers of the $(d-1)$-dimensional hemispheres bounding these lunes belong to $\bd (W)$. 
\end{cla}

The following Theorem from \cite{L5} generalizes a theorem from \cite{LaMu1} for $d=2$. 
The proof presented below is analogous, but this time we apply Theorem \ref{smooth} in the proof. 

\begin{thm} \label{Delta-CW} 
If a reduced convex body $R \subset S^d$ satisfies $\Delta(R) \geq \frac{\pi}{2}$, then $R$ is a body of constant width $\Delta (R)$. 
\end{thm}

Let us return for a while to Example \ref{ExRedBodyOn$S^2$}.
Having in mind the above theorem we observe that if $|ab'| = \frac{\pi}{2}$, then $P$ is a body of constant width $\frac{\pi}{2}$.
In particular, we may take $|ac| = \frac{\pi}{6}$ and $|cb'| = \frac{\pi}{3}$. 
Clearly, if $|ca| = |cb'|$, then our $P$ is a disk.
If $a=c=a'$, then $P$ is a quarter of a disk. 

Recall that every spherical convex body of constant width is reduced. 
Theorem \ref{Delta-CW} says that the opposite implication holds for bodies of thickness at least $\frac{\pi}{2}$, i.e. that every reduced spherical convex body of thickness at least $\frac{\pi}{2}$ is of constant width.
This does not hold for reduced spherical convex bodies of thickness below $\frac{\pi}{2}$ since every regular spherical triangle of thickness below $\frac{\pi}{2}$ is reduced but not of constant width. 

Here is a claim from \cite{L5}. 

\begin{cla} \label{diameterBEITRAGE} 
If for a convex body $C \subset S^d$ the numbers $\Delta(C)$ and $\diam (C)$ are equal and at most $\frac{\pi}{2}$, then $C$ is of constant width $\Delta(C)= \diam (C)$.
\end{cla}

The following two theorems are from \cite{LaMu2}.
 
\begin{thm}\label{touching ball}  
At every boundary point $p$ of a body $W \subset S^d$ of constant width $w > \pi/2$ we can inscribe a unique ball $B_{w- \pi/2}(p')$ touching $W$ from inside at $p$. 
Moreover, $p'$ belongs to the arc connecting $p$ with the center of the unique hemisphere supporting $W$ at $p$, and $|pp'|=w-\frac{\pi}{2}$.
\end{thm}

\begin{thm} \label{strictly} 
Every spherical convex body of constant width smaller than $\frac{\pi}{2}$ on $S^d$ is strictly convex. 
\end{thm}

This theorem was earlier shown in \cite{LaMu1} for $d=2$ only.
Let us add that the thesis of this theorem does not hold for spherical convex bodies of constant width at least $\frac{\pi}{2}$, as we conclude from the following example given in \cite{LaMu1}.
 
\begin{exa} \label{prolong}
{\rm Take a spherical regular triangle $abc$ of sides of length $\kappa < \frac{\pi}{2}$ and prolong them by the same distance $\sigma \leq \frac{\pi}{2} - \kappa$ in both ``directions" up to points $d, e, f, g, h, i$, so that $i, a, b, f$ are on a great circle in this order, $e, b, c, h$ are on a great circle in this order, and $g, c, a, d$ are on a great circle in this order. 
Consider three pieces of circles of radius $\kappa + \sigma$: with center $a$ from $f$ to $g$, with center $b$ from $h$ to $i$, with center $c$ from $d$ to $e$. 
Consider also three pieces of circles of radius $\sigma$: with center $a$ from $i$ to $d$, with center $b$ from $e$ to $f$, with center $c$ from $g$ to $h$. 
The convex hull $U$ of these six pieces of circles is a spherical body of constant width $\kappa + 2\sigma$. 
In particular, when $\kappa + \sigma = \frac{\pi}{2}$, three boundary circles of $U$ become arcs; namely $de$, $fg$ and $hi$.}
\end{exa} 

The following theorem from \cite{LaMu2} gives a positive answer to Problem  \ref{lune passes} in the case of spherical bodies of constant width. 
It is a generalization of the version for $S^2$ given as Theorem 5.3 in \cite{LaMu1}. 
The idea of the proof of our theorem below for $S^d$ substantially differs from the one mentioned for $S^2$.

\begin{thm} \label{center} 
For every body $W \subset S^d$ of constant width $w$ and every $p \in \bd (W)$ there exists a lune $L \supset W$ satisfying $\Delta (L) = w$ with $p$ as the center of one of the two $(d-1)$-dimensional hemispheres bounding this lune.
\end{thm}

If the body $W$ from Theorem \ref{center} is of constant width greater than $\frac{\pi}{2}$, then  by Theorem \ref{smooth} it is smooth. 
Thus at every $p \in \bd (W)$ there is a unique supporting hemisphere of $W$, and so the lune $L$ from the formulation of this theorem is unique.
If the constant width of $W$ is at most $\frac{\pi}{2}$, there are non-smooth bodies of constant width (e.g., a Reuleaux triangle on $S^2$) and then for non-smooth $p \in \bd (W)$ we may have more such lunes. 

Theorem \ref{center} implies the first statement of the following corollary.
The second statement follows from Theorem \ref{smooth} and the last part of Lemma \ref{DistInLunes} from~\cite{L3}.

\begin{cor} \label{pq} 
For every convex body $W \subset S^d$ of constant width $w$ and for every $p \in \bd (W)$ there exists $q \in \bd (W)$ such that $|pq| = w$.  
If $w > \frac{\pi}{2}$, this $q$ is unique.
\end{cor} 

\begin{thm} \label{diam=w} 
If $W \subset S^d$ is a body of constant width $w$, then $\diam (W)=w$.
\end{thm}

The above theorem is proved in \cite{LaMu2} and the following one in \cite{L3}.

\begin{thm} \label{smooth reduced} 
Every smooth reduced body on $S^d$ is of constant width.
\end{thm}

Having in mind that the definition of constant width by Ara\'ujo is equivalent to our definition of constant width (see the last Section), we may present the case of his Theorem B from \cite{Ara} for $S^d$ as follows.

\begin{thm} \label{Araujo} 
If a body $W \subset S^2$ of constant width $w$ has perimeter $p$ and area $a$, then $p = (2\pi - a) \tan \frac{w}{2}$.
\end{thm}

The next proposition proved in \cite{L5} is applied in the proof of the forthcoming Theorem \ref{wulff}.

\begin{pro} \label{equivalent} 
The following conditions are equivalent:

{\rm (1)} \ $C \subset S^d$ is a reduced body with $\Delta(C) = \frac{\pi}{2}$,

{\rm (2)} \ $C \subset S^d$ is a reduced body with $\diam (C) = \frac{\pi}{2}$,

{\rm (3)} \ $C \subset S^d$ is a body of constant width $\frac{\pi}{2}$,

{\rm (4)} \ $C \subset S^d$ is of constant diameter $\frac{\pi}{2}$.
\end{pro}

By part (b) of Theorem \ref{iff}, we conclude the following variant of Theorem \ref{Delta-CW} presented in \cite{L5}.

\begin{cor} \label{ThenConstWidth} 
If a reduced convex body $R \subset S^d$ satisfies $\, \diam(R) \geq \frac{\pi}{2}$, then $R$ is a body of constant width $w$ equal to $\diam (R)$. 
It is also a body of constant diameter $w$.
\end{cor}

Let us comment on an application of spherical bodies of constant width in the research on the Wulff shape.
Recall that Wulff \cite{Wu} defined a geometric model of a crystal equilibrium, later named {\bf Wulff shape}\index{Wulff shape}.  
The literature concerning this and related subjects is very comprehensive. 
For a survey see the monograph \cite{PV} by Pimpinelli and Vilain, also the article \cite {HN2} by Han and Nishimura.  
In \cite{HN1} and \cite{HN3} the authors consider the dual Wulff shape and the {\bf self-dual Wulff shape}\index{self-dual shape} as a Wulff shape which is equal to its dual.
They apply the classical notion of central projection from the open hemisphere centered at a point $n \in S^d$ into the hyperplane of $E^{d+1}$ supporting $S^d$ at $n$.
This hyperplane may be treated as $E^d$ with origin $n$.
The image of a Wulff shape under the inverse projection is called the {\bf spherical convex body induced by this Wulff shape}.
The paper \cite{HN1} presents a proof that a Wulff shape is self-dual if and only if the spherical convex body induced by it is a spherical body of constant width $\frac{\pi}{2}$.
Hence from Proposition \ref{equivalent} we obtain the following theorem.
In particular, the equivalence with the third condition gives a positive answer to a question by Han and Nishimura considered at the end of their paper \cite{HN2}.

\begin{thm} \label{wulff}
{\it Each of the following conditions is equivalent to the statement that the Wulff shape $W_\gamma$ is self-dual:

\noindent
- the spherical convex body induced by $W_\gamma$ is of constant width $\frac{\pi}{2}$,

\noindent
- the spherical convex body induced by $W_\gamma$ is a reduced body of thickness $\frac{\pi}{2}$,

\noindent
- the spherical convex body induced by $W_\gamma$ is a reduced body of diameter $\frac{\pi}{2}$,

\noindent
- the spherical convex body induced by $W_\gamma$ is a body of constant diameter $\frac{\pi}{2}$.}
\end{thm}

Some more applications of the papers \cite{L3}, \cite{L6}, \cite{LaMu2} and \cite{Mus} on the width of spherical convex bodies to Wulff shape are in the paper \cite{Han}.

The last theorem in this section, being a spherical version of the theorem of Blaschke \cite{Bla2}, is a special case of Theorem \ref{approxReduced} for reduced bodies.

\begin{thm} 
Let $W \subset S^2$ be a body of constant width $w \leq \frac{\pi}{2}$. 
For arbitrary $\varepsilon >0$ there exists a body $W_\varepsilon \subset S^2$ of constant width $w$ whose boundary consists only of circles of radius $w$ such that the Hausdorff distance between $W_\varepsilon$ and $W$ is at most $\varepsilon$.
\end{thm}

\vskip0.21cm

\section{Complete spherical convex bodies} \label{Complete}

\vskip0.21cm
\noindent
Similarly to the traditional notion of a complete set in the Euclidean space $E^d$ (for instance, see \cite{BF}, \cite{ChGr},
 \cite{Eg} and \cite{MMO}) we say that any subset of a hemisphere of $S^d$ which is a largest (in the sense of inclusion) set of a given diameter $\delta \in (0, \frac{\pi}{2})$ is a {\bf complete set of diameter $\delta$}, or for brevity, a {\bf complete set}\index{complete set}. 
By the way, the above definition adds the leaking assumption that our set ``is a subset of a hemisphere" to the definition given in \cite{L7}. 
This correction is also added in the arXiv version mentioned in \cite{L7}.

\begin{thm} \label{SubsetOfComplete}
An arbitrary set of diameter $\delta \in (0, \pi)$ on the sphere $S^d$ is a subset of a complete set of diameter $\delta$ on $S^d$.
\end{thm}

The proof of this theorem from  \cite{L7} is similar to the proof by Lebesgue \cite{Leb} for $E^d$ (it is recalled in Part 64 of \cite{BF}).
Let us add that earlier P\'al \cite{Pal} proved this for $E^2$ by a different method.

The following fact from \cite{L7} permits to use the term {\bf complete convex body}\index{complete convex body} for a complete set.

\begin{cla}\label{intersectionBALLS}  
Any complete set of diameter $\delta$ on $S^d$ coincides with the intersection of all balls of radius $\delta$ centered at points of this set.
As a consequence, every spherical complete set is a convex body.
\end{cla}

Below we see two lemmas from  \cite{L7} needed for the proof of Theorem \ref{mainComplete}.

\begin{lem}\label{p'} 
If $C \subset S^d$ is a complete body of diameter $\delta$, then for every $p \in \bd (C)$ there exists $p' \in C$ such that $|pp'|=\delta$.
\end{lem}

For distinct points $a, b \in S^d$ at distance $\delta < \pi$ from a point $c \in S^d$ we define the following piece $P_c(a,b)$ of a circle as the set of points $v \in S^d$ such that $cv$ has length $\delta$ and intersects $ab$.
See Figure \ref{fig:circle} with $P_c(a,b)$ drawn by a broken line.
Also the second lemma, presented below, is formulated for $S^d$ despite the fact that we apply it in the proof of the forthcoming Theorem \ref{mainComplete} only for the case of $S^2$. 

\begin{lem}\label{circle} 
Let $C \subset S^d$ be a complete convex body of diameter $\delta$.
Take $P_c(a,b)$ with $|ac|$ and $|bc|$ equal to $\delta$ such that $a, b \in C$ and $c \in S^d$.
Then $P_c(a,b) \subset C$.
\end{lem}

\begin{figure}[htbp]
\hskip2.5cm \includegraphics[width=8.11cm,height=5.08cm]{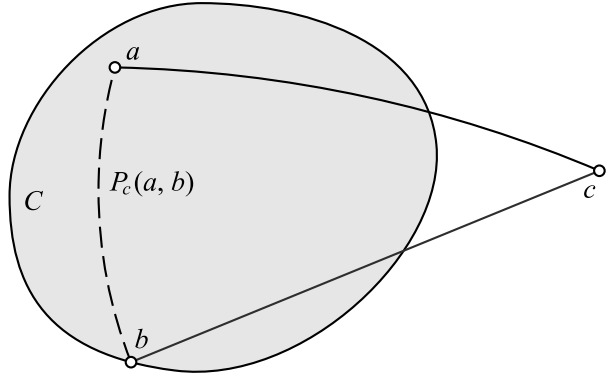} \\  
\caption{Illustration to Lemma \ref{circle}} \label{fig:circle}
\end{figure}

The following theorem from \cite{L7} presents the spherical version of the classical theorem in $E^d$ proved by Meissner \cite {Me} for $d=2,3$ and by Jessen \cite{Je} for arbitrary $d$.

\begin{thm}\label{mainComplete} 
A convex body of diameter $\delta$ on $S^d$ is complete if and only if it is of constant width~$\delta$.
\end{thm}

\begin{proof}
In the first part of the proof let us prove that if a body $C \subset S^d$ of diameter $\delta$ is complete, then it is of constant width $\delta$.

Suppose the contrary, i.e., that $\width_J(C) \not = \delta$ for a hemisphere $J$ supporting $C$.
Our aim is to get a contradiction.

By Theorem  \ref{MaxWidth=Diam} and Proposition \ref{Prop1L2} we have $\width_J(C) \leq \delta$, and so by the fact that $\width_J(C) \not = \delta$ we obtain $\Delta (C) < \delta$.
Recall that $\Delta (C)$ is the thickness of each narrowest lune containing $C$.
Take such a lune $L = G \cap H$. 
Denote by $g, h$ the centers of $G/H$ and $H/G$, respectively (see Figure \ref{fig:complete}).
Of course, $|gh| < \delta$.
By Claim \ref{BothTheCenters} we have $g, h \in C$.
By Lemma \ref{p'} there exists a point $g' \in C$ at distance $\delta$ from $g$.
If $\Delta(L) \geq \frac{\pi}{2}$, then by $C \subset L$ we get $g' \in L$, which contradicts with the last part of Claim \ref{DistInLunes}.

\vskip0.2cm
\begin{figure}[htbp]
\hskip1.05cm \includegraphics[width=10.8cm,height=7.83cm]{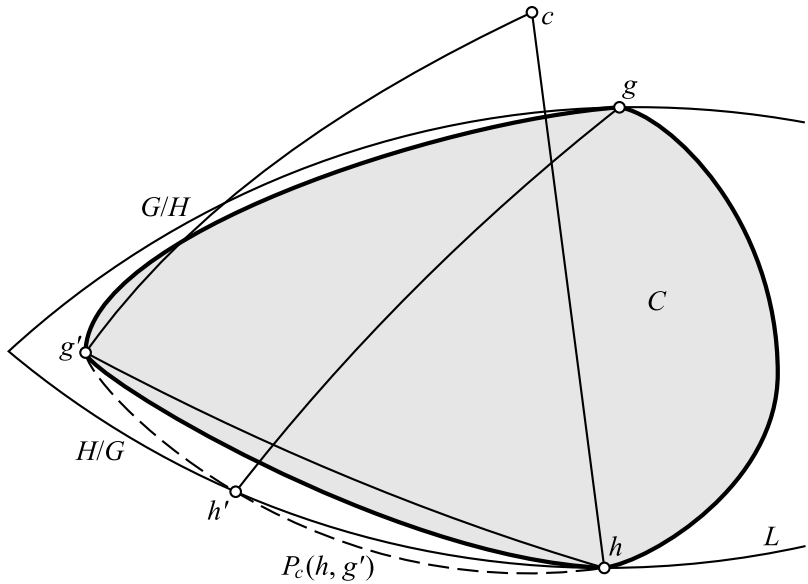} \\ 
\caption{First illustration to the proof of Theorem \ref{mainComplete}} \label{fig:complete}
\end{figure}

Now assume that $\Delta(L) < \frac{\pi}{2}$.
Since the triangle $ghg'$ is non-degenerate, there is a unique two-dimensional sphere $S^2 \subset S^d$ containing $g, h, g'$.
Clearly, $ghg'$ is a subset of $M= C \cap S^2$.
Hence $M$ is a convex body on $S^2$.
Denote by $F$ this hemisphere of $S^2$ such that $hg' \subset \bd(F)$ and $g \in F$.
There is a unique $c \in F$ such that $|ch| = \delta = |cg'|$.
By Lemma \ref{circle} for $d=2$ we have $P_c(h,g') \subset M$ (see Figure \ref{fig:complete}).

We intend to show that $c$ is not on the great circle $J$ of $S^2$ through $g$ and $h$.
In order to see this, suppose for a while the opposite, i.e. that $c \in J$.
Then from $|g'g| = \delta$, $|g'c| = \delta$ and $|hc| = \delta$ we conclude  that $\angle gg'c = \angle hcg'$.
So the spherical triangle $g'gc$ is isosceles, which together with $|gg'|= \delta$ gives $|cg| = \delta$.
Since $|gh| = \Delta(L) = \Delta(C) >0$ and $g$ is a point of $ch$ different from $c$, 
which is impossible.
Hence $c \not \in J$.

By the preceding paragraph $P_c(h,g')$ intersects $\bd(M)$ at a point $h'$ different from $h$ and $g'$.
So the open piece of $P_c(h,g')$ between $h$ and $h'$ is out of $L$ and thus also out of $M$ (see again Figure \ref{fig:complete}).
This gives a contradiction with the result of the paragraph before the last.
Consequently, $C$ is a body of constant width $\delta$.

\vskip0.1cm 
In the second part of the proof let us prove that if $C$ is a spherical body of constant width $\delta$, then $C$ is a complete body of diameter $\delta$.

In order to prove this, it is sufficient to take any point $r \not \in C$ and to show that $\diam (C \cup \{r\}) > \delta$.

\vskip0.2cm
\begin{figure}[htbp]
\hskip2.01cm \includegraphics[width=8.94cm,height=6.44cm]{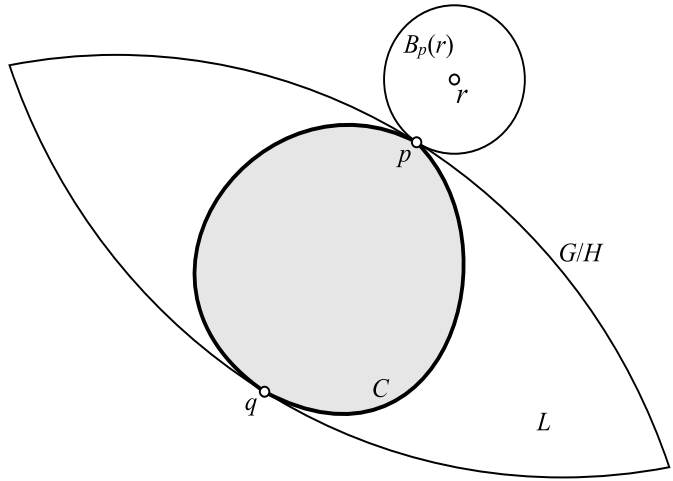} \\ 
\vskip-0.15cm
\caption{Second illustration to the proof of Theorem \ref{mainComplete}} \label{fig:complete+}
\end{figure}

Take the largest ball $B_{p}(r)$ disjoint from the interior of $C$ (see Figure \ref{fig:complete+}).
Since $C$ is convex, $B_{p}(r)$ has exactly one point $p$ in common with $C$.  
By Theorem \ref {center} there exists a lune $L = G \cap H$ containing $C$ of thickness $\delta$ with $p$ as the center of $G/H$.
Denote by $g$ the center of $H/G$. 
By Claim \ref{BothTheCenters} we also have $q \in C$.
From the fact that $rp$ and $pq$ are orthogonal to $G/H$ at $p$, we see that $p \in rq$.
Moreover, $p$ is not an endpoint of $rq$ and $|pq| = \delta$,
Hence $|rq| > \delta$.
Thus $\diam (C \cup \{r\}) > \delta$.
Since $r \not \in C$ is arbitrary, $C$ is complete.
\end{proof}

The above theorem is applied by Bezdek \cite{Bez2} for a generalization of the theorem of Leichtweiss \cite{Lei} saying that  amongst all bodies of constant width $w$ on $S^2$ only the Reauleaux triangle has the smallest possible area, as conjecture by Blaschke \cite{Bla2}).

Here is next theorem from \cite{L7}.

\begin{thm}\label{constant diameter} 
Bodies of constant diameter on $S^d$ coincide with complete bodies.
\end{thm}

\begin{proof}
Take a complete body $D \subset S^d$ of diameter $\delta$.
Let $g \in \bd (D)$ and $G$ be a hemisphere supporting $D$ at $g$.
By Theorem \ref{mainComplete} the body $D$ is of constant width $\delta$.
So $\width_G(D) = \delta$ and there exists a hemisphere $H$ such that the lune $G \cap H \supset D$ has thickness $\delta$.
By Claim \ref{BothTheCenters} the centers $h$ of $H/G$ and $g$ of $G/H$ belong to $D$.
So $|gh| =\delta$.
Thus $D$ is of constant diameter $\delta$.

Consider a body $D \subset S^d$ of constant diameter $\delta$.
Let $r \not \in D$.
Take the largest $B_{\rho}(r)$ whose interior is disjoint from $D$.
Denote by $p$ the common point of $B_{\rho}(r)$ and $D$.
A unique hemisphere $J$ supports $B_{\rho}(r)$ at $p$.
Observe that $D \subset J$ (if not, there is a point $v \in D$ outside $J$; clearly $vp$ passes through $\inter (B_{\rho}(r))$, a contradiction).
Since $D$ is of constant diameter $\delta$, there is $p' \in D$ with $|pp'| =\delta$.
Observe that $\angle rpp' \geq \frac{\pi}{2}$.
If it is equal to $\frac{\pi}{2}$, then $|rp'| > \delta$.
If it is larger than $\frac{\pi}{2}$, the triangle $rpp'$ is obtuse and then by the law of cosines $|rp'| > |pp'|$ and hence
$|rp'| > \delta$.
Since $|rp'| > \delta$ in both cases we see that $D$ is complete.
\end{proof}

Theorem \ref{mainComplete} permits to change ``complete" into ``constant width" in Theorem \ref{constant diameter}, so we get the following one proved in \cite{L7}.
It is analogous to the result in $E^d$ given by Reidemeister \cite{Re}.

\begin{thm} \label{CD=CW}
Bodies of constant diameter on $S^d$ coincide with bodies of constant width.
\end{thm}

Independently, this fact is also proved by Han and Wu in their forthcoming paper \cite{HanWu} providing a quite different consideration using the polarity on $S^d$.
Let us add that earlier same partial facts were established. 
Namely, in \cite{LaMu2} it is shown that any body of constant width $\delta$ on $S^d$ is of constant diameter $\delta$ and the inverse is argued for $\delta \geq \frac{\pi}{2}$, and in \cite{L6} for $\delta < \frac{\pi}{2}$ if $d=2$.

\vskip0.1cm
Only here after presenting a number of facts, we are ready to discuss our latest Problem.
It is a spherical analog of a conjecture from \cite{L1} for reduced bodies in $E^2$. 
There such a conjecture for $E^d$ with $d \geq 3$ is not true by the example on p. 558 in \cite{L0} of a reduced body of thickness $1$ and arbitrarily large fixed diameter.

\begin{prob}
Is it true that every reduced body $R \subset S^d$, where $\Delta (R) < \frac{\pi}{2}$, is a subset of a ball of radius $\Delta(R)$ centered at a boundary point of $R$? 
\end{prob}

The answer to this problem is positive for each body of constant width on $S^d$, what is more for every its boundary point.
This follows by our Theorem \ref{mainComplete} (it follows by Claim \ref{intersectionBALLS} as well).
So by Theorem \ref{Delta-CW} this is also true for every reduced $R \subset S^d$ with $\Delta (R) \geq \frac{\pi}{2}$.
By the way, {Micha\l} Musielak (the author of \cite{Mus} and coauthor of \cite{LaMu1}-\cite{LaMu3}) has got an unpublished yet positive answer to this question for reduced polygons on $S^2$.

\vskip0.21cm

\section{Final remarks} \label{Final}

\vskip0.11cm
Convex spherical sets are considered by many authors.
Mostly they go with the same natural notion of convexity as in this chapter. 
Other definitions are considered very seldom in the literature.
Three soft variations of this natural notion are recalled in Section 9.1 of the article \cite{DGK} by Danzer, Gr\"unbaum and Klee.
The definitions recalled there of a spherically convex set are slightly weaker than the definition of convex sets on $S^d$ used in the present chapter (and named strongly convex sets in \cite{DGK}), i.e., the obtained families of convex sets in these three senses are ``slightly" wider than the family of convex sets resulting from the above natural definition.

We know the following two earlier approaches to define the width of a spherically convex body.

Santal\'o \cite{Sa1} (see also \cite{Sa2}) defines a notion of breadth of a convex body $C \subset S^2$ corresponding to a boundary point of it.
We surmise that he tacitly assumes that $C$ is a smooth body because his definition is as follows.
He takes a boundary point $p$ of $C$, the supporting unique hemisphere $K$ of $C$ at $p$ and the great circle $P$ orthogonal to $K$ at $p$.
Next he finds the other hemisphere $H$, which is orthogonal to $P$ and supports $C$.
By the breadth of $C$ corresponding to $p$, Santal\'o understands the distance from $p$ to the great circle bounding $H$.
Even if we mean this definition as the width generated by the unique hemisphere $K$ supporting $C$, this notion differs from the notion $\width_K(C)$ considered by us; just look to Figure \ref{fig:I-III}.
Gallego, Revent\'os, Solanes, and Teufel \cite{GRST} continue the reserach on the same notion calling it width in place of breadth.

Analogously, Hernandez Cifre and Martinez Fern\'andez \cite{HC-MF} define a notion of width of $C \subset S^2$ at a smooth point $p$ of $\bd (C)$. 
Next they add that if $p$ is an acute point of the boundary of $C$, so if there is a family of more than one supporting hemispheres of $C$ at $p$, then the width of $C$ at $p$ is meant to be the minimum of the above widths for all hemispheres $H$ supporting $C$ at $p$. 
This improved notion still does not have the advantage to be continuous as $p$ moves on the boundary of $C$.
A simple example for this is any regular triangle on $S^2$ whose sides are of length less than $\frac{\pi}{2}$.
Clearly, this notion, as well as the preceding one by S\'antalo, can be defined also for $S^d$.

Recall that our definition of a body of constant width $C \subset S^d$ in Section \ref{Width} is based on the notion of the width $\width_K(C)$ of $C$ determined by a supporting hemisphere $K$ of $C$. 
Surprisingly, some authors define bodies of constant width on the sphere without using any notion of width as the base.
Namely, each of them takes into account a property of a convex body in $E^d$ equivalent of the definition of a body of constant width. 
Such a property is used as the condition defining a convex body on $S^d$ (sometimes only for $d=2$).
Let us recall some such definitions.

In 1915 Sant\'alo \cite{Sa1}, who defined the breadth of spherically convex body at a smooth boundary point (see three paragraphs before), using this notion defines a convex body of constant breadth. 
Namely, he says that a convex body $W \subset S^2$ is of constant breadth if for every boundary point its minimum breadth corresponding to this point and the maximum breadth corresponding to the same point are equal. 
Since he tacitly assumes that $W$ is smooth, below we also make this assumption. 
Clearly, this definition is equivalent to the statement that at every boundary point of $W$ the breadth is the same. 
If a smooth convex body $W$ is of constant width in our sense, then by Claim \ref{determines hemisphere} it is of constant breadth in the sense of Sant\'alo. 
On the other hand, if $W$ is of constant breadth $w$ in the sense of Sant\'alo, then for the unique (by the smoothness) supporting hemisphere $G$ at any fixed boundary point there is a hemisphere $H$ supporting $W$ such that the lune $G \cap H$ has width $w$. 
We know that $H$ touches $\bd (W)$, but we do not know if the midpoints of $G/H$ and $H/G$ are in $W$. 
So we do not know the answer if every body of constant breadth in the sense of Sant\'alo is of constant width in our sense.
Gallego, Revent\'os, Solanes and Teufel \cite{GRST} consider the same notion (which the call constant width) for any smooth convex body in a Riemannian space, so in particular for $S^d$. 

In 1996 Ara\'ujo \cite{Ara} presents his own proof of the planar case of the theorem of Meissner \cite {Me} that bodies of constant diameter coincide with bodies of constant width.
Next, by analogy he defines bodies of constant width on $S^2$ just as bodies of constant diameter.
By our Theorem \ref{constant diameter} his notion is equivalent to our notion of constant width.

In 1983 Chakerian and Groemer \cite{ChGr} calls a convex body $W \subset S^2$ to be of constant width $w$ if it can be rotated in a lune $L$ of thickness $w$ such that $W$ always touches both the bounding semicircles of $L$ at it midpoints (they need this notion in order to recall a result of Blaschke from \cite{Bla1} on a spherical analog of the Barbier's theorem, who rather intuitively tells on constant width on $S^2$). 
In order to avoid the term ``rotate", let us rephrase this condition into ``for every hemisphere $G$ supporting $W$ there exists a hemisphere $H$ such that the lune $G \cap H$ is an orthogonally supporting lune of $W$".
By Claim \ref{lunesCW} we conclude that every body of constant width in our sense is also of constant width in the sense of Chakerian and Groemer. 
Is the converse true?
Taking into account the above rephrased condition, we see that the answer is positive if the answer to the Problem \ref{opposite} is positive.

In 1995 Dekster \cite{Dek1} calls a convex body $C \subset S^d$ of constant width $w>0$ if for any $p \in \bd(W)$ and any normal $m$ to $C$ at $p$ (by such $m$ Dekster means the orthogonal direction to a particular supporting hemisphere of $C$ at $p$) there exists an arc $pq$ of length $w$ having the direction $-m$ such that $pq\subset C$ but $C$ does not contain any longer arc $pq' \supset pq$. 
If $W$ is a body of constant width $w$ in our sense, then by Claim \ref{determines hemisphere} it is of constant width $w$ in the sense of Dekster.
Take a body of constant width $W$ in the sense of Dekster.
Assume that $\diam(W) \leq \frac{\pi}{2}$.
His Theorem 2 from \cite{Dek3} implies that at the ends of the above arc $pq$ the body $W$ is supported by two hemispheres orthogonal to $pq$, which by Claim \ref{C-subset-L} and Theorem \ref{MaxWidth=Diam} leads to the conclusion that $W$ is of constant width in our sense.
To sum up, the families of bodies of constant width in the sense of Dekster and in our sense are identical.

In 2005 Leichtweiss \cite{Lei} considers a so called strip region $T_G$ of width $b$ on $S^2$ determined by a great circle $G$ as the set consisting of all points distant by at most $\frac{b}{2}$ from $G$. 
If a convex body $C$ is contained in a $T_G$ and touched by the two circles bounding it, he says that $T_G$ is a supporting strip region of $C$.
He calls a convex body $W \subset S^2$ of constant width $b$ if all supporting strip regions $T_G$ of $W$ have width $b$.
If a body $W$ is of constant width $w$ in our sense, then from Claim \ref{determines hemisphere} it follows that $W$ is the intersection of the family ${\mathcal L}(W)$ of all lunes supporting $W$, and that the midpoints of the semicircles bounding every lune from ${\mathcal L}(W)$ belong to $W$.
Of course all these lunes have thickness $w$.
Take into account all great circles $G$ of $S^2$ passing through the corners of every lune from $\mathcal L(W)$.
Clearly all these strip regions $T_G$ have the same width $w$.
Consequently, $W$ is of constant width in the sense of Leichtweiss.
Conversely, assume that $W$ is of constant width in the sense of Leichtweiss. 
By virtue of the convexity of $W$, any strip region $T_G$ of $W$ has exactly one point on each of the two bounding circles of $T_G$, and so their distance is $w$.
Hence adding any extra point to $W$ increases the diameter.
Thus $W$ is complete, and by Theorem \ref{mainComplete} it is of constant diameter $w$.
We conclude that the family of bodies of constant width in the sense of Leichtweiss coincides with the family of bodies of constant width in our sense.
Observe that all this can be generalized to $S^d$.
Namely, in place of $G$ take a $(d-1)$-dimensional great subsphere of $S^d$ and as the strip region $T(G)$ of width $b$, again take the set of all points distant by at most $\frac{b}{2}$ from $G$. 

In 2012 Bezdek \cite{Bez1} calls a convex body $W \subset S^d$ of constant width if the smallest thickness of a lune containing $W$ is equal to $\diam (W)$ (see also his recent preprint \cite{Bez2}).
If $W$ is of constant width $w$ in the sense of this chapter, then it is of constant width in the sense of Bezdek.
Still the narrowest lune containing $W$ has thickness $w$, and $\diam(W) = w$ by Theorem \ref{diam=w}.
On the other hand, let us show that if a convex body $W \subset S^d$ is of constant width $w$ in the sense of Bezdek, then $W$ is of constant width in our sense.
Assume the contrary, that is that $W$ is not of constant width in our sense.
Then not all widths of $W$ are equal.
Thus there are hemispheres $H$ and $H'$ supporting $W$ such that $\width_H(W) < \width_{H'}(W)$. 
By Theorem \ref{MaxWidth=Diam} and Proposition \ref{Prop1L2} we have $\width_{H'}(W) \leq \diam(W)$.
So $\width_H(W) < \diam(W)$.
This contradicts Bezdek's condition that the narrowest lune containing $W$ has thickness $\diam(W)$.
Both implications show that the families of constant width in Bezdek's sense and in our sense coincide.

The literature concerning the subject of this paper is very wide. 
Here we chronologically list a few valuable papers on this subject not quoted earlier in this chapter: \cite{Eu}, \cite{Had1}, \cite{Had2}, \cite{Dek2}, \cite{GHS}, \cite{Bak}, \cite{FIN}, \cite{BeS}, \cite{MaMa}, \cite{BeLa}, \cite{BHPS}, \cite{Pap3}, \cite{Guo}, \cite{ZhouGuo}, \cite{GuoPeng} .

\vskip0.2cm

\noindent
Marek Lassak

\vskip-0.1cm
\noindent
University of Technology and Life Sciences 

\vskip-0.1cm
\noindent
al. Kaliskiego 7, 85-796 Bydgoszcz, Poland

\vskip-0.1cm
\noindent
e-mail address: lassak@utp.edu.pl

\end{document}